\documentclass[a4paper]{amsart}
\usepackage{amscd, amssymb, color, booktabs, chemarrow}
\usepackage{myown}

\usepackage{harpoon}

\title[Theta Ranks]{On Theta Ranks of Irreducible Characters of a Finite Classical Group}
\author{Shu-Yen Pan}
\address{Department of Mathematics,
National Tsing Hua University, Hsinchu 300, Taiwan}
\email{sypan@math.nthu.edu.tw}


\date{\today}

\begin{document}

\begin{abstract}
In this paper, we interpret the theta rank of an irreducible character of a finite classical
group in terms of the data from the Lusztig classification.
Then we prove the following two results: (1) the agreement of the $U$-rank and the theta rank
for any irreducible character of low Theta rank; (2) the agreement of the parabolic asymptotic rank
and the theta rank for any irreducible character.
\end{abstract}

\maketitle
\tableofcontents

\section{Introduction}

\subsection{}
Let $\bfG=\bfG_n$ be a finite symplectic group $\Sp_{2n}$,
a finite orthogonal group $\rmO^+_{2n}$, $\rmO^-_{2n+2}$, $\rmO_{2n+1}$,
or a finite unitary group $\rmU_{2n}$, $\rmU_{2n+1}$ defined over a finite field
$\bff_q$ where $q$ is a power of an odd prime.
Let $G$ denote the finite group of rational points of $\bfG$.
For an irreducible character $\rho$ of $G_n$,
several notions of ``rank'' attached to $\rho$ are defined in \cite{gurevich-howe-rank},
namely the \emph{$U$-rank} (denoted by $U\text{\rm -rk}(\rho)$),
the \emph{asymptotic rank} (denoted by $A\text{\rm -rk}(\rho)$),
and the \emph{tensor rank} (denoted by $\otimes\text{\rm -rk}(\rho)$).
It is known that $0\leq U\text{\rm -rk}(\rho)\leq n$.
An irreducible character $\rho$ of $G_n$ is said to be \emph{of low $U$-rank} if
$U\text{\rm -rk}(\rho)<n$ (resp.~$U\text{\rm -rk}(\rho)<n-1$)
if $G_n$ is a symplectic group or a unitary group (resp.~an orthogonal group).

Because the tensor rank of $\rho$ can be interpreted as the first occurrence in the theta correspondence
for some dual pair involving $G$, in this article the tensor rank will be called the
\emph{$\Theta$-rank} and denoted by $\Theta\text{\rm -rk}(\rho)$ (\cf.~Subsection~\ref{0318}).
Note that $\Theta\text{\rm -rk}(\rho)$ is always even if $G$ is an orthogonal group.
A similar notion of the $\Theta$-rank which is called \emph{level} and denoted by
$\grl(\rho)$ is also defined in \cite{GLT01} and \cite{GLT02}.
Due to different normalizations, we can see that
\[
\grl(\rho)=\begin{cases}
\Theta\text{\rm -rk}(\rho), & \text{if $G$ is a symplectic group or a unitary group;}\\
\frac{1}{2}\Theta\text{\rm -rk}(\rho), & \text{if $G$ is an orthogonal group.}
\end{cases}
\]
Moreover, we will modify the definition of the asymptotic ranks and define
the \emph{parabolic asymptotic rank} (denoted by $\overline A\text{\rm -rk}(\rho)$, \cf.~Subsection~\ref{0414}).

Of course, it is expected these ranks are closely related.
In fact, one can see that
\[
U\text{\rm -rk}(\rho)\leq A\text{\rm -rk}(\rho)\leq \overline A\text{\rm -rk}(\rho)
\leq\otimes\text{\rm -rk}(\rho)=\Theta\text{\rm -rk}(\rho)
\]
(\cf.~(\ref{0413})) for any $\rho\in\cale(G)$.
Moreover, Gurevich-Howe make the following conjecture in \cite{gurevich-howe-rank}:
\begin{conj}[Gurevich-Howe]\label{0105}
Let $\rho$ be an irreducible character of a finite classical group $G$.
Then
\begin{enumerate}
\item[(i)] $A\text{\rm -rk}(\rho)=\Theta\text{\rm -rk}(\rho)$;

\item[(ii)] if $\rho$ is of low $U$-rank,
then $U\text{\rm -rk}(\rho)=\Theta\text{\rm -rk}(\rho)$.
\end{enumerate}
\end{conj}
As indicated in \cite{gurevich-howe-rank}, if the conjecture is true,
the theta correspondence provides a method to construct an irreducible character of a
given $U$-rank or $A$-rank.
When $G$ is a general linear group, (i) of the conjecture is already proved in \cite{gurevich-howe-rank} theorem 11.4.
One of the purposes of the article is to provide similar results for other classical groups.

\subsection{}
Let $\cale(G)$ denote the set of irreducible characters of $G$.
It is known that the set $\cale(G)$ partitions into subsets $\cale(G)_s$, called \emph{Lusztig series},
indexed by the set of conjugacy classes of semisimple elements $s$ in the connected component
of the dual group $G^*$ of $G$ (\cf.~\cite{lusztig-book}).
Elements in $\cale(G)_1$ are called \emph{unipotent}.
It is known by Lusztig that unipotent characters $\rho=\rho_\Lambda\in\cale(G)_1$
are parametrized by the equivalence classes of \emph{symbols}
\[
\Lambda=\binom{a_1,a_2,\ldots,a_{m_1}}{b_1,b_2,\ldots,b_{m_2}}\in\cals_\bfG
\]
satisfying certain conditions where $\{a_i\}$, $\{b_j\}$ are finite sets of non-negative integers written
in decreasing order, i.e, $a_1>a_2>\cdots>a_{m_1}$ and $b_1>b_2>\cdots>b_{m_2}$ (\cf.~\cite{lg}).
Note that by replacing $\Lambda$ by another symbol in its equivalence class if necessary,
we can always assume that both rows $\{a_i\},\{b_j\}$ of $\Lambda$ are non-empty.
Our first result (\cf.~Proposition~\ref{0306}) is to interpret the $\Theta$-rank of
a unipotent character $\rho_\Lambda$ in terms of the invariants of symbols:
\begin{equation}\label{0104}
\Theta\text{\rm -rk}(\rho_\Lambda)=\begin{cases}
\min\{2\,{\rm rk}(\Lambda^\rmt\smallsetminus\textstyle\binom{-}{a_1}),
2\,{\rm rk}(\Lambda^\rmt\smallsetminus\textstyle\binom{b_1}{-})\}, & \text{if $\bfG$ is $\Sp$ or $\rmO$};\\
\min\{{\rm rk}_\rmU(\Lambda^\rmt\smallsetminus\textstyle\binom{-}{a_1}),
{\rm rk}_\rmU(\Lambda^\rmt\smallsetminus\textstyle\binom{b_1}{-})\}& \text{if $\bfG$ is $\rmU$}
\end{cases}
\end{equation}
where
\begin{align*}
\binom{a_1,a_2,\ldots,a_{m_1}}{b_1,b_2,\ldots,b_{m_2}}^\rmt
&=\binom{b_1,b_2,\ldots,b_{m_2}}{a_1,a_2,\ldots,a_{m_1}} \\
{\rm rk}(\Lambda)
&=\sum_{i=1}^{m_1}a_i+\sum_{j=1}^{m_2}b_j-\left\lfloor\left(\frac{m_1+m_2-1}{2}\right)^2\right\rfloor,\\
{\rm rk}_\rmU(\Lambda)
&=\sum_{i=1}^{m_1} 2a_i+\sum_{j=1}^{m_2} 2b_j
+\frac{|m_1-m_2|}{2}-\frac{(m_1+m_2)(m_1+m_2-2)}{2}.
\end{align*}
From (\ref{0104}), we see that $\Theta\text{\rm -rk}(\rho_\Lambda)$ is always even if $\bfG$ is a symplectic
group or an orthogonal group.

\subsection{}
Next we want to investigate the $\Theta$-rank of a general irreducible character.
First suppose $\bfG$ is a symplectic group or an orthogonal group.
For $s\in(G^*)^0$, we can define $\bfG^{(0)},\bfG^{(-)},\bfG^{(+)}$ (\cf.~Subsection~\ref{0215}) where
$\bfG^{(0)}$ is a product of unitary groups or general linear groups,
and $\bfG^{(-)},\bfG^{(+)}$ are given as follows:
\begin{itemize}
\item if $\bfG=\rmO^\epsilon_{2n}$,
then $\bfG^{(-)}=\rmO^{\epsilon'}_{2n^{(-)}}$, and $\bfG^{(+)}=\rmO^{\epsilon'\epsilon}_{2n^{(+)}}$;

\item if $\bfG=\Sp_{2n}$,
then $\bfG^{(-)}=\rmO^{\epsilon'}_{2n^{(-)}}$, and $\bfG^{(+)}=\Sp_{2n^{(+)}}$;

\item if $\bfG=\rmO_{2n+1}$,
then $\bfG^{(-)}=\Sp_{2n^{(-)}}$, and $\bfG^{(+)}=\Sp_{2n^{(+)}}$,
\end{itemize}
for some $\epsilon'=\pm$, and $n^{(-)}, n^{(+)}$ depending on $s$.
Then we have a (modified) \emph{Lusztig correspondence}
\[
\Xi_s\colon \cale(G)_s\rightarrow \cale(G^{(0)}\times G^{(-)}\times G^{(+)})_1.
\]
The mapping $\Xi_s$ is two-to-one: $\Xi_s(\rho)=\Xi_s(\rho\cdot\sgn)$ when $\bfG$ is an odd-orthogonal group;
and $\Xi_s$ is bijective otherwise.
Then we can write
\[
\Xi_s(\rho)=\rho^{(0)}\otimes\rho_{\Lambda^{(-)}}\otimes\rho_{\Lambda^{(+)}}
\]
for $\rho\in\cale(G)_s$ where $\rho^{(0)}\in\cale(G^{(0)})_1$, $\Lambda^{(-)}\in\cals_{\bfG^{(-)}}$,
and $\Lambda^{(+)}\in\cals_{\bfG^{(+)}}$.

Next suppose $\bfG$ is a unitary group.
As described in Subsection~\ref{0316},
for the Lusztig correspondence $\grL_s\colon\cale(G)_1\rightarrow\cale(C_{G^*}(s))_1$
we can write
\begin{equation}
\grL_s(\rho)=\bigotimes_{i=1}^r\bigotimes_{j=1}^{t_i}\rho^{(ij)}.
\end{equation}
Note that $t_1=q+1$ and $\rho^{(1,j)}$ is a unipotent character of $\rmU_{n_{1,j}}(q)$
for $j=1,\ldots,q+1$.
Now we have the following formula (\cf.~Proposition~\ref{0307} and Proposition~\ref{0311})
for the $\Theta$-rank of a general irreducible character in terms of data from
the Lusztig classification:
\begin{thm}\label{0103}
Let $\rho\in\cale(G)$.
\begin{enumerate}
\item[(i)] If\/ $\bfG=\Sp_{2n}$,
then
\[
\Theta\text{\rm -rk}(\rho)
=\min\{2n-2n^{(+)}+\Theta\text{\rm -rk}(\rho_{\Lambda^{(+)}}),2n-2n^{(-)}+\Theta\text{\rm -rk}(\rho_{\Lambda^{(-)}})+1\}.
\]
\item[(ii)] If\/ $\bfG=\rmO^\epsilon_{2n}$ or $\rmO_{2n+1}$,
then
\[
\Theta\text{\rm -rk}(\rho)
=\min\{2n-2n^{(+)}+\Theta\text{\rm -rk}(\rho_{\Lambda^{(+)}}),2n-2n^{(-)}+\Theta\text{\rm -rk}(\rho_{\Lambda^{(-)}})\}.
\]
\item[(iii)] If\/ $\bfG=\rmU_n$,
then
\[
\Theta\text{\rm -rk}(\rho)
=\min\{n-n_{1,j}+\Theta\text{\rm -rk}(\rho^{(1,j)})\mid j=1,\ldots,q+1\}.
\]
\end{enumerate}
\end{thm}

\subsection{}
After we analyze several properties of the $\Theta$-ranks,
the following result can be derived easily from the result on $\eta$-correspondence for a reductive dual
pair in stable range by Gurevich-Howe in \cite{gurevich-howe-rank}:
\begin{thm}\label{0101}
Let $\rho\in\cale(G_n)$ where $\bfG_n$ is $\Sp_{2n}$, $\rmO^+_{2n}$, $\rmO^-_{2n+2}$, $\rmO_{2n+1}$,
$\rmU_{2n}$, or $\rmU_{2n+1}$.
If\/ $\Theta\text{\rm -rk}(\rho)\leq n$,
then
\[
U\text{\rm -rk}(\rho)
=A\text{\rm -rk}(\rho)
=\overline A\text{\rm -rk}(\rho)
=\Theta\text{\rm -rk}(\rho).
\]
\end{thm}

Finally, by using the same arguments in \cite{gurevich-howe-rank} section 11,
we can show the agreement of the parabolic asymptotic rank and the theta rank:

\begin{thm}\label{0102}
Let $\rho\in\cale(G)$ where $G$ is a symplectic group, an orthogonal group or a unitary group.
Then
\[
\overline A\text{\rm -rk}(\rho)=\Theta\text{\rm -rk}(\rho).
\]
\end{thm}
Theorems~\ref{0101} and \ref{0102} provide a partial answer for Conjecture~\ref{0105}.

\subsection{}
The contents of this article are as follows.
In Section~2 we first review some results on the Lusztig's classification on unipotent characters
and the Lusztig correspondence.
Then we prove several results on the first occurrence for theta correspondence.
In Section~3, we define the $\Theta$-rank of an irreducible character via the theta correspondence.
Then we investigate several properties of the $\Theta$-ranks and prove Theorem~\ref{0103}.
We discuss the relation of $\Theta$-ranks and other ranks in the final section.
In particular, Theorem~\ref{0101} and Theorem~\ref{0102} are proved.

\section{Theta Correspondence and Lusztig Correspondence}

\subsection{Basic notations}
Let $\bff_q$ be a finite field of $q$ elements where $q$ is a power of an odd prime $p$,
and let $\overline\bff_q$ denote a fixed algebraic closure of $\bff_q$.

Let $\bfG$ be a symplectic group, an orthogonal group or a unitary group defined over $\bff_q$,
and let $G$ denote the finite group of rational points.
The character of an irreducible representation of $G$ is called an \emph{irreducible character}
of $G$, and the set of all irreducible characters of $G$ is denoted by $\cale(G)$.
The character of a one-dimensional representation is also called a \emph{linear character}.
It is proved by Lusztig that the set $\cale(G)$ admits a partition
\[
\cale(G)=\bigcup_{(s)\subset(G^*)^0}\cale(G)_s
\]
indexed by conjugacy classes $(s)$ of semisimple simple elements $s$ in the
connected component $(G^*)^0$ of the dual group $G^*$ of $G$.
Each $\cale(G)_s$ is called a \emph{Lusztig series},
and the elements in $\cale(G)_1$ are called \emph{unipotent}.
The trivial character of $G$ is denoted by $\bf1$.

\subsection{Unipotent characters of a symplectic group or an orthogonal group}
The set of \emph{partitions} of $n$ is denoted by $\calp(n)$.
If $\mu\in\calp(n)$, we also write $|\mu|=n$.
An ordered pair $\sqbinom{\mu}{\nu}$ of two partitions $\mu,\nu$,
is called a \emph{bi-partition of $n$} if $\left|\sqbinom{\mu}{\nu}\right|:=|\mu|+|\nu|=n$.
The set of bi-partitions of $n$ is denoted by $\calp_2(n)$.

Now we recall the notations of symbols from \cite{lg}.
Note that for our convenience, the definition here is slightly different from the original one.
A \emph{$\beta$-set} $A=\{a_1,a_2,\ldots,a_m\}$ is a finite subset (possibly empty) of non-negative integers
written in strictly decreasing order, i.e., $a_1>a_2>\cdots>a_m$.
A \emph{symbol}
\[
\Lambda=\binom{A}{B}=\binom{a_1,a_2,\ldots,a_{m_1}}{b_1,b_2,\ldots,b_{m_2}}
\]
is an ordered pair of two $\beta$-sets $A,B$.
On the set of symbols we define an equivalence relation ``$\sim$'' generated by
\[
\binom{a_1,a_2,\ldots,a_{m_1}}{b_1,b_2,\ldots,b_{m_2}}\sim
\binom{a_1+1,a_2+1\ldots,a_{m_1}+1,0}{b_1+1,b_2+1,\ldots,b_{m_2}+1,0}.
\]
The set of equivalence classes of symbols is denoted by $\cals$.

For a symbol $\Lambda=\binom{a_1,\ldots,a_{m_1}}{b_1,\ldots,b_{m_2}}$,
we define its \emph{defect} and \emph{rank} by
\begin{align*}
{\rm def}(\Lambda) &=m_1-m_2, \\
{\rm rk}(\Lambda) &=\sum_{i=1}^{m_1}a_i+\sum_{j=1}^{m_2}b_j
-\left\lfloor\left(\frac{m_1+m_2-1}{2}\right)^2\right\rfloor.
\end{align*}
It is not difficult to see that two equivalent symbols have the same defect and the same rank.
For a symbol $\Lambda=\binom{A}{B}$ we define its \emph{transpose} $\Lambda^\rmt=\binom{B}{A}$.
It is obvious that ${\rm rk}(\Lambda^\rmt)={\rm rk}(\Lambda)$ and
${\rm def}(\Lambda^\rmt)=-{\rm def}(\Lambda)$.
We define a mapping $\Upsilon$ from symbols to bi-partitions by
\begin{equation}
\Upsilon\colon\binom{a_1,\ldots,a_{m_1}}{b_1,\ldots,b_{m_2}}\mapsto
\sqbinom{a_1-(m_1-1),a_2-(m_1-2),\ldots,a_{m_1-1}-1,a_{m_1}}
{b_1-(m_2-1),b_2-(m_2-2),\ldots,b_{m_2-1}-1,b_{m_2}}.
\end{equation}
The image $\Upsilon(\Lambda)$ in fact depends only on the equivalence class of $\Lambda$.
Moreover, it is easy to check that
\[
{\rm rk}(\Lambda)=\begin{cases}
|\Upsilon(\Lambda)|+\tfrac{1}{4}({\rm def}(\Lambda)-1)({\rm def}(\Lambda)+1),
& \text{if ${\rm def}(\Lambda)$ is odd};\\
|\Upsilon(\Lambda)|+\tfrac{1}{4}({\rm def}(\Lambda))^2,
& \text{if ${\rm def}(\Lambda)$ is even}.
\end{cases}
\]

For a symplectic group or an orthogonal group $\bfG$,
we define the following sets of (equivalence classes of) symbols:
\begin{align*}
\cals_{\rmO^+_{2n}}
&= \{\,\Lambda\in\cals\mid{\rm rk}(\Lambda)=n,\ {\rm def}(\Lambda)\equiv 0\pmod 4\,\}; \\
\cals_{\Sp_{2n}}
&= \{\,\Lambda\in\cals\mid{\rm rk}(\Lambda)=n,\ {\rm def}(\Lambda)\equiv 1\pmod 4\,\}; \\
\cals_{\rmO^-_{2n}}
&= \{\,\Lambda\in\cals\mid{\rm rk}(\Lambda)=n,\ {\rm def}(\Lambda)\equiv 2\pmod 4\,\}; \\
\cals_{\SO_{2n+1}}=\cals_{\rmO_{2n+1}}
&= \{\,\Lambda\in\cals\mid{\rm rk}(\Lambda)=n,\ {\rm def}(\Lambda)\equiv 3\pmod 4\,\}.
\end{align*}
The following are modified from the results proved by Lusztig (\cf.~\cite{lg}):
\begin{enumerate}
\item Suppose that $\bfG=\Sp_{2n},\rmO^\epsilon_{2n},\SO_{2n+1}$ where $\epsilon=+$ or $-$.
There is a parametrization $\cale(G)_1$ by $\cals_\bfG$,
so a unipotent character $\rho$ of $G$ is written as $\rho=\rho_\Lambda$ for $\Lambda\in\cals_\bfG$.
If $\bfG=\rmO^\epsilon_{2n}$, it is known that $\rho_\Lambda\cdot\sgn=\rho_{\Lambda^\rmt}$.

\item Suppose that $\bfG=\rmO_{2n+1}$.
Then a unipotent character $\rho\in\cale(G)_1$ can be written as $\rho=\rho_\Lambda$ or $\rho=\rho_\Lambda\cdot\sgn$
for $\Lambda\in\cals_{\SO_{2n+1}}$
\end{enumerate}

\subsection{Unipotent characters of a unitary group}
In this subsection, suppose that $\bfG$ is a unitary group.
For a symbol $\Lambda$, we define
\begin{equation}\label{0208}
{\rm rk}_\rmU(\Lambda)=2|\Upsilon(\Lambda)|+\tfrac{1}{2}|{\rm def}(\Lambda)|(|{\rm def}(\Lambda)|+1).
\end{equation}
It is easy to check that
\begin{equation}
{\rm rk}_\rmU\binom{a_1,a_2,\ldots,a_{m_1}}{b_1,b_2,\ldots,b_{m_2}}
=\sum_{i=1}^{m_1} 2a_i+\sum_{j=1}^{m_2} 2b_j
+\frac{|m_1-m_2|}{2}-\frac{(m_1+m_2)(m_1+m_2-2)}{2}.
\end{equation}

Then we let $\cals_{\rmU_n}$ be the set of (equivalence classes of) symbols $\Lambda\in\cals$ such that
\begin{itemize}
\item ${\rm def}(\Lambda)$ is either even and non-negative, or odd and negative;

\item ${\rm rk}_\rmU(\Lambda)=n$.
\end{itemize}
Again, there is a parametrization $\cale(\rmU_n(q))_1$ by $\cals_{\rmU_n}$ (\cf.~\cite{FS} or \cite{pan-eta-unitary}),
i.e., for $\rho\in\cale(\rmU_n(q))_1$, we can write $\rho=\rho_\Lambda$ for a unique $\Lambda\in\cals_{\rmU_n}$.

\subsection{Lusztig correspondence}\label{0215}
The following result is proved by Lusztig (\cf.~\cite{DM} theorem 13.23, remark 13.24):

\begin{prop}[Lusztig]
There is a bijection
\[
\grL_s\colon\cale(G)_s\rightarrow\cale(C_{G^*}(s))_1
\]
such that $R^\bfG_{\bfT^*,s}\mapsto \varepsilon_\bfG\varepsilon_{C_{\bfG^*}(s)}R^{C_{\bfG^*}(s)}_{\bfT^*,1}$
where $C_{G^*}(s)$ denotes the centralizer in $G^*$ of $s$,
$\varepsilon_\bfG=(-1)^\kappa$ and $\kappa$ is the relative rank of\/ $\bfG$.
\end{prop}

The bijection $\grL_s$ called a \emph{Lusztig correspondence}.

\begin{rem}
Clearly, ${\rm def}(\Lambda)\equiv 1\pmod 4$ if and only if ${\rm def}(\Lambda^\rmt)\equiv 3\pmod 4$.
The the Lusztig correspondence $\grL_1\colon\cale(\Sp_{2n}(q))_1\rightarrow\cale(\SO_{2n+1}(q))_1$ is a
bijection given by $\rho_\Lambda\mapsto\rho_{\Lambda^\rmt}$.
\end{rem}

For $s\in (G^*)^0$, we define
\begin{align}\label{0203}
\begin{split}
G^{(0)}=G^{(0)}(s)
&=\prod_{\langle\lambda\rangle\subset\{\lambda_1,\ldots,\lambda_n\},\ \lambda\neq\pm 1}G_{[\lambda]}(s);\\
G^{(-)}=G^{(-)}(s)
&=G_{[-1]}(s); \\
G^{(+)}=G^{(+)}(s)
&=\begin{cases}
(G_{[1]}(s))^*, & \text{if $G$ is symplectic};\\
G_{[1]}(s), & \text{otherwise},
\end{cases}
\end{split}
\end{align}
where $G_{[\lambda]}(s)$ is given in \cite{amr} subsection 1.B
(see also \cite{pan-Lusztig-correspondence} subsection 2.2).
Hence we have a bijection
\[
\cale(C_{G^*}(s))_1\simeq\begin{cases}
\cale(G^{(0)}\times G^{(-)}\times G^{(+)})_1\times\{\pm1\}, & \text{if $G$ is odd-orthogonal};\\
\cale(G^{(0)}\times G^{(-)}\times G^{(+)})_1, & \text{otherwise}.
\end{cases}
\]
It is known that $\bfG^{(0)}$ is a product of unitary groups or general linear groups,
and $\bfG^{(-)},\bfG^{(+)}$ are given as follows:
\begin{itemize}
\item if $\bfG=\rmO^\epsilon_{2n}$,
then
$\bfG^{(-)}=\rmO^{\epsilon'}_{2n^{(-)}}$, and $\bfG^{(+)}=\rmO^{\epsilon'\epsilon}_{2n^{(+)}}$
for some $\epsilon'=+$ or $-$;

\item if $\bfG=\Sp_{2n}$,
then $\bfG^{(-)}=\rmO^{\epsilon'}_{2n^{(-)}}$, and $\bfG^{(+)}=\Sp_{2n^{(+)}}$
for some $\epsilon'=+$ or $-$;

\item if $\bfG=\rmO_{2n+1}$,
then $\bfG^{(-)}=\Sp_{2n^{(-)}}$, and $\bfG^{(+)}=\Sp_{2n^{(+)}}$;

\item if $\bfG=\rmU_{n}$,
then $\bfG^{(-)}=\rmU_{n^{(-)}}$, and $\bfG^{(+)}=\rmU_{n^{(+)}}$,
\end{itemize}
where $n^{(+)}=n^{(+)}(s)$ (resp.~$n^{(-)}=n^{(-)}(s)$) denotes the multiplicity
of $+1$ (resp.~$-1$) in the set $\{\lambda_1,\ldots,\lambda_n\}$ of ``eigenvalues''
of $s$.

So now a Lusztig correspondence $\grL_s$ induces a mapping
\[
\Xi_s\colon \cale(G)_s\rightarrow \cale(G^{(0)}(s)\times G^{(-)}(s)\times G^{(+)}(s))_1
\]
called the \emph{modified Lusztig correspondence}.
The mapping $\Xi_s$ is two-to-one and $\Xi_s(\rho)=\Xi_s(\rho\cdot\sgn)$ when $\bfG$ is an odd-orthogonal group;
and $\Xi_s$ is a bijection otherwise.
We can write $\Xi_s(\rho)=\rho^{(0)}\otimes\rho^{(-)}\otimes\rho^{(+)}$.
Because now $\rho^{(\varepsilon)}$ where $\varepsilon=+$ or $-$ is a unipotent character of a symplectic group,
an even-orthogonal group, or a unitary group,
we can write $\rho^{(\varepsilon)}=\rho_{\Lambda^{(\varepsilon)}}$ for some symbol
$\Lambda^{(\varepsilon)}\in\cals_{\bfG^{(\varepsilon)}}$.
Now $\rho\mapsto\Lambda^{(\varepsilon)}$ is a mapping from $\cale(G)_s$ to $\cals_{\bfG^{(\varepsilon)}}$
and $\Lambda^{(\varepsilon)}$ is also denoted by $\Lambda^{(\varepsilon)}(\rho)$,
and hence we can write
\[
\Xi_s(\rho)=\rho^{(0)}\otimes\rho_{\Lambda^{(-)}}\otimes\rho_{\Lambda^{(+)}}.
\]
Note that when $\bfG$ is an orthogonal group,
$\rho\in\cale(G)_s$ if and only if $\rho\cdot\sgn\in\cale(G)_s$.

\begin{lem}
Let $G$ be an orthogonal group and $\rho\in\cale(G)_s$ for some $s$.
Then
\[
\Lambda^{(\varepsilon)}(\rho\cdot\sgn)=
\begin{cases}
\Lambda^{(\varepsilon)}(\rho)^\rmt, & \text{if $G$ is an even orthogonal group};\\
\Lambda^{(\varepsilon)}(\rho), & \text{if $G$ is an odd orthogonal group}.
\end{cases}
\]
\end{lem}
\begin{proof}
This is obvious from the definition in (\ref{0203}).
\end{proof}

Recall that a linear character $\chi_\bfG$ of order two is defined when
$\bfG$ is an orthogonal group or a unitary group (\cf.~\cite{pan-odd} subsection~3.3).
Then there exists a bijection $\cale(G)_s\rightarrow\cale(G)_{-s}$ given by
$\rho\mapsto\rho\chi_\bfG$.

\begin{lem}\label{0212}
Let $G$ be an orthogonal group or a unitary group, and let $s$ be a semisimple element in $(G^*)^0$.
Then $G^{(0)}(-s)\simeq G^{(0)}(s)$, $G^{(-)}(-s)\simeq G^{(+)}(s)$ and $G^{(+)}(-s)\simeq G^{(-)}(s)$.
\end{lem}
\begin{proof}
This is obvious from the definition in (\ref{0203}).
\end{proof}

\begin{lem}\label{0207}
Let $G$ be an orthogonal group, and let $\rho\in\cale(G)$.
Then $\Lambda^{(\varepsilon)}(\rho\chi_\bfG)=\Lambda^{(-\varepsilon)}(\rho)$
for $\varepsilon=+$ or $-$.
\end{lem}
\begin{proof}
Let $\rho\in\cale(G)_s$ for some semisimple $s\in (G^*)^0$.
It is known that $\chi_\bfG R^{\bfG}_{\bfT^*,s}=R^{\bfG}_{\bfT^*,-s}$ for any rational
maximal torus $\bfT^*$ in $\bfG^*$ containing $s$.
Then from Lemma~\ref{0212} we have $\bfG^{(-)}(-s)\simeq\bfG^{(+)}(s)$ and $\bfG^{(+)}(-s)\simeq\bfG^{(-)}(s)$.
If we write $\Xi_s(\rho)=\rho^{(0)}\otimes\rho^{(-)}\otimes\rho^{(+)}$ and
$\Xi_{-s}(\rho\chi_\bfG)=\rho'^{(0)}\otimes\rho'^{(-)}\otimes\rho'^{(+)}$,
then we have $\rho'^{(-)}=\rho^{(+)}$ and $\rho'^{(+)}=\rho^{(-)}$.
This means that $\Lambda^{(-)}(\rho\chi_\bfG)=\Lambda^{(+)}(\rho)$
and $\Lambda^{(+)}(\rho\chi_\bfG)=\Lambda^{(-)}(\rho)$.
\end{proof}

\subsection{Lusztig correspondence and parabolic induction}\label{0211}
Let $\bfG_n$ denote the group of split rank $n$ in its Witt series,
i.e., $\bfG_n$ is one of $\Sp_{2n}$, $\rmO^+_{2n}$, $\rmO^-_{2n+2}$, $\rmO_{2n+1}$,
$\rmU_{2n}$, or $\rmU_{2n+1}$.
For $m\geq n$,
let $\cale^{G_m}_{G_n}(\rho)$ denote the set of irreducible components of $\Ind_{P_n}^{G_m}(\rho\otimes\bf1)$
where $P_n$ is a parabolic subgroup of $G_m$ whose Levi factor is isomorphic to $G_n\times T_{m-n}$
and $T_{m-n}$ is the direct product of $m-n$ copies of $\GL_1(q)$ if $\bfG_n$
is a symplectic group or an orthogonal group;
$T_{m-n}$ is the direct product of $m-n$ copies of $\GL_1(q^2)$ if $\bfG_n$ is a unitary group.
For $m<n$ and $\rho\in\cale(G_n)$,
let $\cale_{G_n}^{G_m}(\rho)$ denote the set of $\rho'\in\cale(G_m)$ such that
$\rho\in\cale_{G_m}^{G_n}(\rho')$.

Let $\rho\in\cale(G_n)$ and $\rho_m\in\cale^{G_m}_{G_n}(\rho)$ for some $m\geq n$.
Because
\[
\Ind^{G_m}_{P_n}(R_{\bfT^*,s}^{\bfG_n}\otimes{\bf1})
=R^{\bfG_m}_{\bfT^*\times\bfT_{m-n},s_m}
\]
where $s_m=(s,1)\in (G^*_n)^0\times T_{m-n}\subset (G^*_m)^0$,
we see that $\langle \rho_m,R^{\bfG_m}_{\bfT^*\times\bfT_{m-n},s_m}\rangle\neq 0$
if and only if $\langle \rho,R^{\bfG_n}_{\bfT^*,s}\rangle\neq 0$.
This means that if $\rho\in\cale(G_n)_{s}$ for some $s$,
then $\rho_m\in\cale(G_m)_{s_m}$.
In particular, if $\rho$ is unipotent,
then every element in $\cale_{G_n}^{G_m}(\rho)$ is unipotent.

\begin{lem}\label{0202}
Suppose that $m\geq n$, $\rho\in\cale(G_n)_s$ for some $s$,
and $\Xi_s(\rho)=\rho^{(0)}\otimes\rho^{(-)}\otimes\rho^{(+)}$.
Then
\[
\cale_{G_n}^{G_m}(\rho)
=\{\,\rho_m\in\cale(G_m)_{s_m}\mid\Xi_{s_m}(\rho_m)=\rho^{(0)}\otimes\rho^{(-)}\otimes\rho_m^{(+)},\
\rho^{(+)}_m\in\cale_{G_{n^{(+)}}}^{G_{m^{(+)}}}(\rho^{(+)})\,\}
\]
where $m^{(+)}=n^{(+)}+(m-n)$.
\end{lem}
\begin{proof}
Let $\rho\in\cale(G_n)_s$.
If $\rho_m\in\cale_{G_n}^{G_m}(\rho)$,
then from above we see that $\rho_m\in\cale(G_m)_{s_m}$ where $s_m=(s,1)\in (G^*_n)^0\times T_{m-n}\subset(G_m^*)^0$.
Then $\bfG^{(0)}(s_m)=\bfG^{(0)}(s)$,
$\bfG^{(-)}(s_m)=\bfG^{(-)}(s)$,
and $\bfG^{(+)}(s_m)=\bfG_{m^{(+)}}$ where $m^{(+)}=n^{(+)}+(m-n)$.
It is known that $R^{\bfG_{m}}_{\bfT^*\times\bfT_{m-n},s_m}=\Ind_{P_n}^{G_m}(R^{G_n}_{\bfT^*,s}\otimes{\bf1})$
where $P_n$ is a parabolic subgroup of $G_m$ whose Levi factor is isomorphic to $G_n\times T_{m-n}$.
Then the lemma follows easily.
\end{proof}

Suppose that $\rho\in\cale(G_n)_1$.
We write $\rho=\rho_\Lambda$ ($\rho=\rho_\Lambda$ or $\rho=\rho_\Lambda\cdot\sgn$ if $\bfG_n=\rmO_{2n+1}$)
for $\Lambda\in\cals_{\bfG_n}$.
From the Littlewood-Richardson rule, we know that an element $\rho'\in\cale_{G_n}^{G_{n+1}}(\rho)$
is of the form $\rho_{\Lambda'}$ ($\rho_{\Lambda'}$ or $\rho_{\Lambda'}\cdot\sgn$ if $\bfG_n=\rmO_{2n+1}$)
where $\Lambda'\in\cals_{\bfG_{n+1}}$ such that
${\rm def}(\Lambda')={\rm def}(\Lambda)$ and $\Upsilon(\Lambda')$ is obtained from
$\Upsilon(\Lambda)=\sqbinom{\mu}{\nu}$ by adding a box to some row (including an empty row)
of $\mu$ or $\nu$.

\begin{exam}
Suppose that $\rho\in\cale(\Sp_{2n})_s$ and $\Xi_s(\rho)=\rho^{(0)}\otimes\rho_{\Lambda^{(-)}}\otimes\rho_{\Lambda^{(+)}}$
such that $\rho_{\Lambda^{(+)}}$ is cuspidal.
Hence $n^{(+)}=t(t+1)$ for some non-negative integer $t$ and
\[
\Lambda^{(+)}=\begin{cases}
\binom{t,t-1,\ldots,1,0}{-}, & \text{if $t$ is even};\\
\binom{-}{t,t-1,\ldots,1,0}, & \text{if $t$ is odd}.
\end{cases}
\]
Hence the parabolic induced character $\Ind_{P_n}^{\Sp_{2(n+1)}}(\rho\otimes{\bf1})$ is the sum of two
irreducible characters
\[
\Xi_{(s,1)}^{-1}(\rho^{(0)}\otimes\rho_{\Lambda^{(-)}}\otimes\rho_{\Lambda_1})
+\Xi_{(s,1)}^{-1}(\rho^{(0)}\otimes\rho_{\Lambda^{(-)}}\otimes\rho_{\Lambda_2})
\]
where $\Xi_{(s,1)}$ is the bijection
\[
\Xi_{(s,1)}\colon\cale(\Sp_{2(n+1)}(q))_{(s,1)}\rightarrow
\cale(G^{(0)}\times G^{(-)}\times\Sp_{2(t(t+1)+1)}(q))_1
\]
and
\[
\{\Lambda_1,\Lambda_2\}=\begin{cases}
\bigl\{\binom{t+1,t-1,t-2,\ldots,1,0}{-},\binom{t+1,t,t-1,\ldots,1,0}{1}\bigr\}, & \text{if $t$ is even};\\
\bigl\{\binom{-}{t+1,t-1,t-2,\ldots,1,0},\binom{1}{t+1,t,t-1,\ldots,1,0}\bigr\}, & \text{if $t$ is odd}.
\end{cases}
\]
\end{exam}

\subsection{Finite Howe correspondences on unipotent characters}\label{0213}
First we recall the concept of finite reductive dual pairs from \cite{howe-finite}.
Let $(\bfG',\bfG)$ be a finite reductive dual pair of one of the following types:
\begin{enumerate}
\item[(I)] $(\rmU_{n'},\rmU_n)$;

\item[(II)] $(\rmO_{2n'}^\epsilon,\Sp_{2n})$ or $(\Sp_{2n'},\rmO^\epsilon_{2n})$ where $\epsilon=+$ or $-$;

\item[(III)] $(\rmO_{2n'+1},\Sp_{2n})$ or $(\Sp_{2n'},\rmO_{2n+1})$
\end{enumerate}
where $n',n$ are non-negative integers.
By restricting the Weil character with respect to a nontrivial additive character $\psi$ of $\bff_q$
to the dual pair $(\bfG',\bfG)$, we obtain a decomposition
\[
\omega_{\bfG',\bfG}^\psi=\sum_{\rho'\in\cale(G'),\ \rho\in\cale(G)}m_{\rho',\rho}\rho'\otimes\rho
\]
where $m_{\rho',\rho}=0,1$.
Then we have a relation between $\cale(G')$ and $\cale(G)$:
\[
\Theta_{\bfG',\bfG}
=\{\,(\rho',\rho)\in\cale(G')\times\cale(G)\mid m_{\rho',\rho}\neq 0\,\},
\]
called the \emph{$\Theta$-correspondence} for the dual pair $(\bfG',\bfG)$.
For $\rho'\in\cale(G')$, we write
\[
\Theta(\rho')
=\Theta_\bfG(\rho')
=\{\,\rho\in\cale(G)\mid(\rho',\rho)\in\Theta_{\bfG',\bfG}\,\}.
\]

It is well known that for cases (I) and (II) above, the unipotent characters
are preserved by $\Theta$-correspondence, i.e., when $(\rho',\rho)\in\Theta_{\bfG',\bfG}$,
then $\rho'\in\cale(G')_1$ if and only if $\rho\in\cale(G)_1$ (\cf.~\cite{adams-moy}).
For these two cases,
the unipotent part $\omega_{\bfG',\bfG,1}$ of the Weil character $\omega_{\bfG',\bfG}^\psi$ has the
decomposition
\[
\omega_{\bfG',\bfG,1}
=\sum_{\rho'\in\cale(G')_1,\ \rho\in\cale(G)_1}m_{\rho',\rho}\rho'\otimes\rho
=\sum_{(\Lambda',\Lambda)\in\calb_{\bfG',\bfG}}\rho_{\Lambda'}\otimes\rho_\Lambda
\]
where the relation $\calb_{\bfG',\bfG}$ are given in \cite{pan-finite-unipotent}
and \cite{pan-Lusztig-correspondence}.

\subsection{Howe correspondence and Lusztig correspondence}\label{0204}
The following result on the commutativity between the (probably chosen) modified Lusztig correspondence
and the Howe correspondence is from \cite{pan-Lusztig-correspondence} theorem 6.9 and theorem 7.9:
\begin{prop}\label{0214}
Let $(\bfG',\bfG)$ be a finite reductive dual pair, and let $\rho'\in\cale(G')_{s'}$,
$\rho\in\cale(G)_s$ for some $s',s$.
Write $\Xi_{s'}(\rho')=\rho'^{(0)}\otimes\rho_{\Lambda'^{(-)}}\otimes\rho_{\Lambda'^{(+)}}$ and
$\Xi_s(\rho)=\rho^{(0)}\otimes\rho_{\Lambda^{(-)}}\otimes\rho_{\Lambda^{(+)}}$.

\begin{enumerate}
\item[(i)] Suppose that\/ $(\bfG',\bfG)=(\rmO^\epsilon_{2n'},\Sp_{2n})$ or $(\rmU_{n'},\rmU_n)$.
Then $(\rho',\rho)$ occurs in $\Theta_{\bfG',\bfG}$ if and only if
\begin{itemize}
\item $\bfG'^{(0)}\simeq\bfG^{(0)}$ and $\rho'^{(0)}=\rho^{(0)}$;

\item $\bfG'^{(-)}\simeq\bfG^{(-)}$ and $\rho_{\Lambda'^{(-)}}=\rho_{\Lambda^{(-)}}$;

\item $(\rho_{\Lambda'^{(+)}},\rho_{\Lambda^{(+)}})$ occurs in $\Theta_{\bfG'^{(+)},\bfG^{(+)}}$,
\end{itemize}
i.e., the diagram
\[
\begin{CD}
\rho' @> \Theta_{\bfG',\bfG} >> \rho \\
@V \Xi_{s'} VV @VV \Xi_s V \\
\rho'^{(0)}\otimes\rho_{\Lambda'^{(-)}}\otimes\rho_{\Lambda'^{(+)}}
@> \id\otimes\id\otimes\Theta_{\bfG'^{(+)},\bfG^{(+)}} >>  \rho^{(0)}\otimes\rho_{\Lambda^{(-)}}\otimes\rho_{\Lambda^{(+)}}
\end{CD}
\]
commutes.

\item[(ii)] Suppose that\/ $(\bfG',\bfG)=(\rmO_{2n'+1},\Sp_{2n})$.
Then $(\rho',\rho)$ occurs in $\Theta_{\bfG',\bfG}$ if and only if
\begin{itemize}
\item $\bfG'^{(0)}\simeq\bfG^{(0)}$ and $\rho'^{(0)}=\rho^{(0)}$;

\item $(\rho_{\Lambda'^{(+)}},\rho_{\Lambda^{(-)}})$ occurs in $\Theta_{\bfG'^{(+)},\bfG^{(-)}}$,

\item $\bfG'^{(-)}\simeq\bfG^{(+)}$ and $\rho_{\Lambda'^{(-)}}=\rho_{\Lambda^{(+)}}$;
\end{itemize}
i.e., the diagram
\[
\begin{CD}
\rho' @> \Theta_{\bfG',\bfG} >> \rho \\
@V \iota\circ\Xi_{s'} VV @VV \Xi_s V \\
\rho'^{(0)}\otimes\rho_{\Lambda'^{(+)}}\otimes\rho_{\Lambda'^{(-)}}
@> \id\otimes\Theta_{\bfG'^{(+)},\bfG^{(-)}}\otimes\id >> \rho^{(0)}\otimes\rho_{\Lambda^{(-)}}\otimes\rho_{\Lambda^{(+)}}
\end{CD}
\]
commutes where $\iota(\rho'^{(0)}\otimes\rho_{\Lambda'^{(-)}}\otimes\rho_{\Lambda'^{(+)}})
=\rho'^{(0)}\otimes\rho_{\Lambda'^{(+)}}\otimes\rho_{\Lambda'^{(-)}}$.
\end{enumerate}
\end{prop}

Recall that an irreducible character $\rho$ of $G=\Sp_{2n}(q)$ is called \emph{pseudo-unipotent}
if $\rho\in\cale(G)_s$ such that $C_{\SO_{2n+1}}(s)=\rmO^\epsilon_{2n}$ for some $\epsilon=+$ or $-$.
Now from the proposition we have several lemmas concerning the first occurrences of
unipotent or pseudo-unipotent characters.

\begin{lem}\label{0303}
Let $\rho$ be an irreducible unipotent character of\/ $\Sp_{2n}(q)$.
Suppose that $\rho$ occurs in the $\Theta$-correspondence for some dual pair
$(\rmO_{2n'+1},\Sp_{2n})$.
Then $n'\geq n$.
\end{lem}
\begin{proof}
Consider the dual pair $(\bfG',\bfG)=(\rmO_{2n'+1},\Sp_{2n})$ and the following commutative diagram:
\[
\begin{CD}
\rho' @> \Theta_{\bfG',\bfG} >> \rho \\
@V \iota\circ\Xi_{s'} VV @VV \Xi_s V \\
\rho'^{(0)}\otimes\rho_{\Lambda'^{(+)}}\otimes\rho_{\Lambda'^{(-)}}
@> \id\otimes\Theta_{\bfG'^{(+)},\bfG^{(-)}}\otimes\id >> \rho^{(0)}\otimes\rho_{\Lambda^{(-)}}\otimes\rho_{\Lambda^{(+)}}
\end{CD}
\]
Because $\rho$ is a unipotent character of $\Sp_{2n}(q)$,
we see that $\bfG^{(0)}=\rmU_0$, $\bfG^{(-)}=\rmO_0^+$, and $\bfG^{(+)}=\Sp_{2n}$.
This means by Proposition~\ref{0214} that $\bfG'^{(-)}=\Sp_{2n}$ and hence $n'\geq n$.
\end{proof}

\begin{lem}\label{0209}
Let $\rho$ be an irreducible pseudo-unipotent character of\/ $\Sp_{2n}(q)$.
Suppose that $\rho$ occurs in the $\Theta$-correspondence for some dual pair
$(\rmO^\epsilon_{2n'},\Sp_{2n})$.
Then $n'\geq n$.
\end{lem}
\begin{proof}
Consider the dual pair $(\bfG',\bfG)=(\rmO^\epsilon_{2n'},\Sp_{2n})$ and the following commutative diagram:
\[
\begin{CD}
\rho' @> \Theta_{\bfG',\bfG} >> \rho \\
@V \Xi_{s'} VV @VV \Xi_s V \\
\rho'^{(0)}\otimes\rho_{\Lambda'^{(-)}}\otimes\rho_{\Lambda'^{(+)}}
@> \id\otimes\id\otimes\Theta_{\bfG'^{(+)},\bfG^{(+)}} >>  \rho^{(0)}\otimes\rho_{\Lambda^{(-)}}\otimes\rho_{\Lambda^{(+)}}
\end{CD}
\]
Because $\rho$ is a pseudo-unipotent character of $\Sp_{2n}(q)$,
we see that $\bfG^{(0)}=\rmU_0$, $\bfG^{(-)}=\rmO_{2n}^{\epsilon'}$ for some $\epsilon'=+$ or $-$,
and $\bfG^{(+)}=\Sp_0$.
This means by Proposition~\ref{0214} that $\bfG'^{(-)}=\rmO^{\epsilon'}_{2n}$ and hence $n'\geq n$.
\end{proof}

\begin{lem}\label{0205}
Let $\rho$ be an irreducible unipotent character of\/ $\rmO^\epsilon_{2n}(q)$.
Suppose that $\rho\chi_{\rmO^\epsilon_{2n}}$ occurs in the correspondence for the dual pair
$(\Sp_{2n'},\rmO^\epsilon_{2n})$.
Then $n'\geq n$.
\end{lem}
\begin{proof}
Consider the dual pair $(\bfG',\bfG)=(\Sp_{2n'},\rmO^\epsilon_{2n})$ and the following commutative diagram:
\[
\begin{CD}
\rho' @> \Theta_{\bfG',\bfG} >> \rho\chi_\bfG \\
@V \Xi_{s'} VV @VV \Xi_s V \\
\rho'^{(0)}\otimes\rho_{\Lambda'^{(-)}}\otimes\rho_{\Lambda'^{(+)}}
@> \id\otimes\id\otimes\Theta_{\bfG'^{(+)},\bfG^{(+)}} >>  \rho^{(0)}\otimes\rho_{\Lambda^{(-)}}\otimes\rho_{\Lambda^{(+)}}
\end{CD}
\]
Because $\rho$ is a unipotent character of $\rmO^\epsilon_{2n}(q)$,
now $\rho\chi_\bfG$ is in $\cale(\rmO_{2n}^\epsilon(q))_{-1}$ and
$\bfG^{(0)}=\rmU_0$, $\bfG^{(-)}=\rmO^\epsilon_{2n}$, and $\bfG^{(+)}=\rmO^+_0$.
This means by Proposition~\ref{0214} that $\bfG'^{(-)}=\rmO^\epsilon_{2n}$ and hence $n'\geq n$.
\end{proof}

\begin{lem}\label{0206}
Let $\rho$ be an irreducible unipotent character of\/ $\rmO_{2n+1}(q)$.
Suppose that $\rho\chi_{\rmO_{2n+1}}$ occurs in the correspondence for some dual pair $(\Sp_{2n'},\rmO_{2n+1})$.
Then $n'\geq n$.
\end{lem}
\begin{proof}
Consider the dual pair $(\bfG',\bfG)=(\Sp_{2n'},\rmO_{2n+1})$ and the following commutative diagram:
\[
\begin{CD}
\rho' @> \Theta_{\bfG',\bfG} >> \rho\chi_{\bfG} \\
@V \Xi_{s'} VV @VV \iota\circ\Xi_s V \\
\rho'^{(0)}\otimes\rho_{\Lambda'^{(-)}}\otimes\rho_{\Lambda'^{(+)}}
@> \id\otimes\Theta_{\bfG'^{(-)},\bfG^{(+)}}\otimes\id >>  \rho^{(0)}\otimes\rho_{\Lambda^{(+)}}\otimes\rho_{\Lambda^{(-)}}
\end{CD}
\]
Because $\rho$ is a unipotent character of $\rmO_{2n+1}(q)$,
now $\rho\chi_\bfG$ is in $\cale(\rmO_{2n+1}(q))_{-1}$ and
$\bfG^{(0)}=\rmU_0$, $\bfG^{(-)}=\Sp_{2n}$, and $\bfG^{(+)}=\Sp_0$.
This means by Proposition~\ref{0214} that $\bfG'^{(+)}=\Sp_{2n}$ and hence $n'\geq n$.
\end{proof}

\begin{lem}\label{0210}
Let $\rho$ be an irreducible unipotent character of\/ $\rmU_n(q)$,
and let $\chi$ be a non-trivial linear character of\/ $\rmU_n(q)$.
Suppose that $\rho\chi$ occurs in the correspondence for some dual pair $(\rmU_{n'},\rmU_{n})$.
Then $n'\geq n$.
\end{lem}
\begin{proof}
Consider the dual pair $(\bfG',\bfG)=(\rmU_{n'},\rmU_n)$ and the following commutative diagram:
\[
\begin{CD}
\rho' @> \Theta_{\bfG',\bfG} >> \rho\chi \\
@V \Xi_{s'} VV @VV \Xi_s V \\
\rho'^{(0)}\otimes\rho_{\Lambda'^{(-)}}\otimes\rho_{\Lambda'^{(+)}}
@> \id\otimes\id\otimes\Theta_{\bfG'^{(+)},\bfG^{(+)}} >> \rho^{(0)}\otimes\rho_{\Lambda^{(-)}}\otimes\rho_{\Lambda^{(+)}}
\end{CD}
\]
Now $\rho$ is unipotent and $\chi$ is a nontrivial linear character,
we see that $\rho\chi$ is in $\cale(G)_s$ where $s$ does not have eigenvalue $1$,
i.e., $\bfG^{(+)}=\rmU_0$.
Because now $\bfG'^{(0)}=\bfG^{(0)}$ and $\bfG'^{(-)}=\bfG^{(-)}$ by Proposition~\ref{0214},
we conclude that $n'\geq n$.
\end{proof}

\section{$\Theta$-ranks of Irreducible Characters}

\subsection{Definition of $\Theta$-ranks}\label{0318}
Let $G$ be a symplectic group, an orthogonal group or a unitary group.
For $\rho\in\cale(G)$,
the \emph{$\Theta$-rank} of $\rho$, denoted by $\Theta\text{\rm -rk}(\rho)$,
is defined as follows:
\begin{enumerate}
\item If $G$ is a symplectic group, we consider dual pairs $(\rmO^\epsilon_k,\bfG)$ and define
\[
\Theta\text{\rm -rk}(\rho)
=\min\{\,k\mid\rho\in\Theta(\rho')\text{ for some }\rho'\in\cale(\rmO^\epsilon_k(q)),\ \epsilon=+\text{ or }-\,\}.
\]

\item If $G$ is an orthogonal group, we consider dual pairs $(\Sp_{2k},\bfG)$ and define
\begin{multline*}
\Theta\text{\rm -rk}(\rho)
=\min\{\,2k\mid\rho\chi\in\Theta(\rho')
\text{ for some }\rho'\in\cale(\Sp_{2k}(q)) \\
\text{ and some linear character }\chi\in\cale(G)\,\}.
\end{multline*}
So the $\Theta$-rank of an irreducible character of an orthogonal group is always even.
Note that for a non-trivial orthogonal group,
there are four linear characters: $\bf1$, $\sgn$, $\chi_\bfG$, and $\sgn\chi_\bfG$.

\item If $G$ is a unitary group, we consider dual pairs $(\rmU_k,\bfG)$ and define
\begin{multline*}
\Theta\text{\rm -rk}(\rho)
=\min\{\,k\mid\rho\chi\in\Theta(\rho')\text{ for some }\rho'\in\cale(\rmU_k(q)) \\
\text{and linear character }\chi\in\cale(G)\,\}.
\end{multline*}
\end{enumerate}
From the above definition, we see that
$\Theta\text{\rm -rk}(\rho)=\Theta\text{\rm -rk}(\rho\chi)$
for any linear character $\chi$ of $G$.

\begin{exam}
By convention, we have $\omega_{\bfG',\bfG}^\psi={\bf 1}_{\bfG'}\otimes{\bf 1}_\bfG$
when $(\bfG',\bfG)=(\Sp_0,\rmO^\epsilon_n)$, $(\rmO_0^+,\Sp_{2n})$ or $(\rmU_0,\rmU_n)$.
Therefore linear characters are the only irreducible characters of $G$ of $\Theta$-rank 0.
\end{exam}

\begin{rem}\label{0304}
Let $\rho$ be an irreducible character of $G$.
\begin{enumerate}
\item Suppose that $\bfG=\Sp_{2n}$.
Let $n_0^\epsilon(\rho)$ denote the minimum of $2n'$ such that
$\rho$ occurs in the $\Theta$-correspondence for the dual pair $(\rmO^\epsilon_{2n'},\Sp_{2n})$.
It is known from \cite{pan-Lusztig-correspondence} that
\[
n^+_0(\rho)+n^-_0(\rho)\leq 4n+2.
\]
Because now both $n^+_0(\rho),n^-_0(\rho)$ are even,
we see that at least one of $n^+_0(\rho),n^-_0(\rho)$ is less than or equal to $2n$.
Moreover, $\Theta\text{\rm -rk}(\rho)\leq\min\{n^+_0(\rho),n^-_0(\rho)\}$ by the definition
of $\Theta$-ranks,
so we have $0\leq\Theta\text{\rm -rk}(\rho)\leq 2n$.

\item Suppose that $\bfG=\rmO^\epsilon_{2n}$ or $\rmO_{2n+1}$.
Let $n_0(\rho)$ denote the minimum of $2n'$ such that
$\rho$ occurs in the $\Theta$-correspondence for the dual pair $(\Sp_{2n'},\bfG)$.
It is known that
\[
n_0(\rho)+n_0(\rho\cdot\sgn)\leq\begin{cases}
4n, & \text{if $\bfG=\rmO^\epsilon_{2n}$};\\
4n+2, & \text{if $\bfG=\rmO_{2n+1}$}.
\end{cases}
\]
Now both $n_0(\rho),n_0(\rho\cdot\sgn)$ are even.
It again implies that $0\leq\Theta\text{\rm -rk}(\rho)\leq 2n$.

\item Suppose that $\bfG=\rmU_n$.
Let $n_0^+(\rho)$ (resp.\ $n_0^-(\rho)$) denote the minimum of $2n'$ (resp.\ $2n'+1$)
such that $\rho$ occurs in the $\Theta$-correspondence for the dual pair
$(\rmU_{2n'},\rmU_n)$ (resp.\ $(\rmU_{2n'+1},\rmU_n)$).
It is known that
\[
n^+_0(\rho)+n^-_0(\rho)\leq 2n+1.
\]
This implies that $0\leq\Theta\text{\rm -rk}(\rho)\leq n$.
\end{enumerate}
\end{rem}

\subsection{$\Theta$-ranks of unipotent or pseudo-unipotent characters}

Let $\Lambda=\binom{a_1,a_2,\ldots,a_{m_1}}{b_1,b_2,\ldots,b_{m_2}}\in\cals_\bfG$
where $\bfG$ is a symplectic group, an orthogonal group, or a unitary group.
Because
$\binom{a_1,a_2,\ldots,a_{m_1}}{b_1,b_2,\ldots,b_{m_2}}
\sim\binom{a_1+1,a_2+1,\ldots,a_{m_1}+1,0}{b_1+1,b_2+1,\ldots,b_{m_2}+1,0}$
and ${\rm rk}(\Lambda)={\rm rk}(\Lambda')$, ${\rm rk}_\rmU(\Lambda)={\rm rk}_\rmU(\Lambda')$
if $\Lambda\sim\Lambda'$,
we may assume that both rows $\{a_i\},\{b_i\}$ are non-empty.
Then we define
\begin{equation}\label{0201}
\Theta\text{\rm -rk}(\Lambda)=\begin{cases}
\min\{2{\rm rk}(\Lambda^\rmt\smallsetminus\textstyle\binom{-}{a_1}),
2{\rm rk}(\Lambda^\rmt\smallsetminus\textstyle\binom{b_1}{-})\}, & \text{if $\bfG$ is $\Sp$ or $\rmO$};\\
\min\{{\rm rk}_\rmU(\Lambda^\rmt\smallsetminus\textstyle\binom{-}{a_1}),
{\rm rk}_\rmU(\Lambda^\rmt\smallsetminus\textstyle\binom{b_1}{-})\}& \text{if $\bfG$ is $\rmU$}
\end{cases}
\end{equation}
where $\Lambda^\rmt\smallsetminus\textstyle\binom{-}{a_1}=\binom{b_1,b_2,\ldots,b_{m_2}}{a_2,\ldots,a_{m_1}}$, etc.
Note that $\Theta\text{\rm -rk}(\Lambda)$ depends only on the equivalence class of $\Lambda$.

\begin{prop}\label{0306}
Let $\bfG$ be a symplectic group, an orthogonal group, or a unitary group,
and let $\rho\in\cale(G)_1$.
Write $\rho=\rho_\Lambda$ ($\rho=\rho_\Lambda$ or $\rho=\rho_\Lambda\cdot\sgn$ if\/ $\bfG=\rmO_{2n+1}$) for some symbol $\Lambda\in\cals_\bfG$.
Then
\[
\Theta\text{\rm -rk}(\rho_\Lambda)=\Theta\text{\rm -rk}(\Lambda).
\]
\end{prop}
\begin{proof}
Let $\Lambda=\binom{a_1,a_2,\ldots,a_{m_1}}{b_1,b_2,\ldots,b_{m_2}}\in\cals_{\bfG}$ and $\rho=\rho_\Lambda$.

\begin{enumerate}
\item Suppose that $\bfG=\Sp_{2n}$.
By Lemma~\ref{0303}, to compute $\Theta\text{\rm -rk}(\rho_\Lambda)$
we only need to consider the dual pairs $(\rmO^\epsilon_{2n'},\Sp_{2n})$ for $\epsilon=+$ or $-$,
i.e., we do not need to consider the first occurrence of $\rho_\Lambda$ in the correspondence
for the dual pair $(\rmO_{2n'+1},\Sp_{2n})$.
From \cite{pan-Lusztig-correspondence} section 8,
we know that the minimum of $2n'$ such that
$\rho_\Lambda$ occurs in the $\Theta$-correspondence for the dual pair $(\rmO^+_{2n'},\Sp_{2n})$
is equal to $2{\rm rk}\binom{b_1,b_2,\ldots,b_{m_2}}{a_2,a_3,\ldots,a_{m_1}}$.
Similarly, the minimum of $2n'$ such that $\rho_\Lambda$ occurs in the $\Theta$-correspondence
for the dual pair $(\rmO^-_{2n'},\Sp_{2n})$ is equal to
$2{\rm rk}\binom{b_2,b_3,\ldots,b_{m_2}}{a_1,a_2,\ldots,a_{m_1}}$.

\item Suppose that $\bfG=\rmO^\epsilon_{2n}$.
By Lemma~\ref{0205}, if $\rho_\Lambda\chi_\bfG$ or $\rho_\Lambda\cdot\sgn\chi_\bfG$
occurs in the $\Theta$-correspondence for the dual pair $(\Sp_{2n'},\rmO^\epsilon_{2n})$,
then $n'\geq n$.
So we need only to consider the first occurrences of $\rho_\Lambda$ and
$\rho_{\Lambda^\rmt}=\rho_\Lambda\cdot\sgn$
in the correspondence for the dual pair $(\Sp_{2n'},\rmO^\epsilon_{2n})$.
We know that the minimum of $2n'$ such that $\rho_\Lambda$ occurs in the $\Theta$-correspondence
for the dual pair $(\Sp_{2n'},\rmO^\epsilon_{2n})$ is equal to
$2{\rm rk}\binom{b_1,b_2,\ldots,b_{m_2}}{a_2,a_3,\ldots,a_{m_1}}$ when $\epsilon=+$;
and equal to $2{\rm rk}\binom{b_2,b_3,\ldots,b_{m_2}}{a_1,a_2,\ldots,a_{m_1}}$ when $\epsilon=-$.
Similarly, the minimum of $2n'$ such that $\rho_\Lambda\cdot\sgn$ occurs in the $\Theta$-correspondence
for the dual pair $(\Sp_{2n'},\rmO^\epsilon_{2n})$ is equal to
$2{\rm rk}\binom{b_2,b_3,\ldots,b_{m_2}}{a_1,a_2,\ldots,a_{m_1}}$ when $\epsilon=+$;
and equal to $2{\rm rk}\binom{b_1,b_2,\ldots,b_{m_2}}{a_2,a_3,\ldots,a_{m_1}}$ when $\epsilon=-$.

\item Suppose that $\bfG=\rmO_{2n+1}$.
By Lemma~\ref{0206}, we know that if $\rho_\Lambda\chi_\bfG$ or $\rho_\Lambda\cdot\sgn\chi_\bfG$
occurs in the $\Theta$-correspondence for the dual pair $(\Sp_{2n'},\rmO_{2n+1})$,
then $n'\geq n$.
So we need only to consider the first occurrences of $\rho_\Lambda$ and
$\rho_\Lambda\cdot\sgn$ in the correspondence for the dual pair
$(\Sp_{2n'},\rmO_{2n+1})$.

Now we have a commutative diagram
\[
\begin{CD}
\rho' @> \Theta_{\bfG',\bfG} >> \rho_\Lambda \\
@V \Xi_{s'} VV @VV \iota\circ\Xi_1 V \\
{\bf 1}\otimes\rho_{\Lambda'^{(-)}}\otimes{\bf1} @> \id\otimes\Theta_{\bfG'^{(-)},\bfG^{(+)}}\otimes\id >>
{\bf1}\otimes\rho_{\Lambda^\rmt}\otimes{\bf 1}
\end{CD}
\]
Now $\bfG^{(+)}=\Sp_{2n}$ and $\bfG'^{(-)}=\rmO^+_{2n'}$,
so the minimum of $2n'$ such that $\rho_\Lambda$ occurs in the $\Theta$-correspondence
for the dual pair $(\Sp_{2n'},\rmO_{2n+1})$ is equal to the minimum of $2n''$ such that
$\rho_{\Lambda^\rmt}$ occurs in the $\Theta$-correspondence for the dual pair
$(\rmO^+_{2n''},\Sp_{2n})$, and hence equal to
$2{\rm rk}\binom{b_1,b_2,\ldots,b_{m_2}}{a_2,a_3,\ldots,a_{m_1}}$.

Moreover, $\Xi_1(\rho_\Lambda\cdot\sgn)=\Xi_1(\rho)$ and
now $\bfG^{(+)}=\Sp_{2n}$ and $\bfG'^{(-)}=\rmO^-_{2n'}$,
so the minimum of $2n'$ such that $\rho_\Lambda\cdot\sgn$ occurs in the $\Theta$-correspondence
for the dual pair $(\Sp_{2n'},\rmO_{2n+1})$ is equal to the minimum of $2n''$ such that
$\rho_{\Lambda^\rmt}\in\cale(\Sp_{2n}(q))_1$ occurs in the $\Theta$-correspondence for the dual pair
$(\rmO^-_{2n''},\Sp_{2n})$ and hence equal to
$2{\rm rk}\binom{a_1,a_2,\ldots,a_{m_1}}{b_2,b_3,\ldots,b_{m_2}}$.

\item Suppose that $\bfG=\rmU_n$.
By Lemma~\ref{0210}, to compute $\Theta\text{\rm -rk}(\rho_\Lambda)$,
we do not need to consider the first occurrences of $\rho_\Lambda\chi$
for any nontrivial linear character $\chi$ of $G$.
Now the minimum of $n'$ such that $\rho_\Lambda$ occurs in the $\Theta$-correspondence
for the dual pair $(\rmU_{n'},\rmU_n)$ is equal to
${\rm rk}_\rmU\binom{b_1,b_2,\ldots,b_{m_2}}{a_2,a_3,\ldots,a_{m_1}}$ if $n'+n$ is even;
and is equal to ${\rm rk}_\rmU\binom{b_2,b_3,\ldots,b_{m_2}}{a_1,a_2,\ldots,a_{m_1}}$ if $n'+n$ is odd.
\end{enumerate}
\end{proof}

The following two corollaries follow from the proof of the above proposition immediately.

\begin{cor}\label{0305}
The $\Theta$-rank of an irreducible unipotent character of\/ $\Sp_{2n}(q)$ is always even.
\end{cor}

Note that if $\rho_\Lambda\in\cale(\rmO_{2n+1}(q))$,
then $\rho_{\Lambda^\rmt}\in\cale(\Sp_{2n}(q))$.

\begin{cor}\label{0310}
If $\rho_\Lambda\in\cale(\rmO_{2n+1}(q))_1$,
then
$\Theta\text{\rm -rk}(\rho_\Lambda)
=\Theta\text{\rm -rk}(\rho_\Lambda\cdot\sgn)
=\Theta\text{\rm -rk}(\rho_{\Lambda^\rmt})$.
\end{cor}

\begin{exam}\label{0317}
It is clear from the result in \cite{adams-moy} and Lemma~\ref{0205} that
the $\Theta$-ranks of cuspidal unipotent characters are described as follows:
\begin{enumerate}
\item The $\Theta$-rank of the unique cuspidal unipotent character of $\Sp_{2d(d+1)}(q)$ is $2d^2$.

\item The $\Theta$-ranks of the two cuspidal unipotent characters of $\rmO^\epsilon_{2d^2}(q)$
(for $d\geq 1$) are $2d(d-1)$ where $\epsilon=+$ if $d$ is even, and $\epsilon=-$ if $d$ is odd.

\item The $\Theta$-ranks of the two cuspidal unipotent characters of $\rmO_{2d(d+1)+1}(q)$ are $2d^2$.

\item The $\Theta$-rank of the unique cuspidal unipotent character of $\rmU_{\frac{1}{2}d(d+1)}(q)$ is $\frac{1}{2}d(d-1)$.
\end{enumerate}
\end{exam}

\begin{exam}
Now we want to compute the $\Theta$-rank of a Steinberg character of $G$.
\begin{enumerate}
\item Suppose that $\bfG=\Sp_{2n}$.
The Steinberg character $\St_{\Sp_{2n}}$ is the unipotent character associated to the symbol
$\binom{n,n-1,\ldots,0}{n,n-1,\ldots,1}$.
Therefore,
\[
\Theta\text{\rm -rk}(\St_{\Sp_{2n}})
= \min\{2{\rm rk}\textstyle\binom{n,n-1,\ldots,1}{n-1,n-2,\ldots,0},
2{\rm rk}\binom{n-1,n-2,\ldots,1}{n,n-1,\ldots,0}\}
= 2n,
\]
i.e., $\St_{\Sp_{2n}}$ is a character of maximal $\Theta$-rank.

\item Suppose that $\bfG=\rmO^+_{2n}$ and $n\geq 1$.
The two Steinberg characters $\rho,\rho\cdot\sgn$ are the unipotent characters associated to the symbols
$\binom{n,n-1,\ldots,1}{n-1,n-2,\ldots,0}$ or $\binom{n-1,n-2,\ldots,0}{n,n-1,\ldots,1}$.
Therefore
\[
\Theta\text{\rm -rk}(\rho)
=\Theta\text{\rm -rk}(\rho\cdot\sgn)
=\min\{2{\rm rk}\textstyle\binom{n-1,n-2,\ldots,0}{n-1,n-2,\ldots,1},
2{\rm rk}\binom{n-2,n-3,\ldots,0}{n,n-1,\ldots,1}\}
=2n-2.
\]
Suppose that $\bfG=\rmO^-_{2n}$ and $n\geq 1$.
The Steinberg characters are the unipotent characters associated to the symbols
$\binom{n,n-1,\ldots,0}{n-1,n-2,\ldots,1}$ or $\binom{n-1,n-2,\ldots,1}{n,n-1,\ldots,0}$.
Therefore,
\[
\Theta\text{\rm -rk}(\rho)
=\Theta\text{\rm -rk}(\rho\cdot\sgn)
=\min\{2{\rm rk}\textstyle\binom{n-1,n-2,\ldots,1}{n-1,n-2,\ldots,0},
2{\rm rk}\binom{n-2,n-3,\ldots,1}{n,n-1,\ldots,0}\}
=2n-2.
\]

\item Suppose that $\bfG=\rmO_{2n+1}$.
The two Steinberg characters $\rho,\rho\cdot\sgn$ of $G$ are the two unipotent characters associated
to the symbol $\binom{n,n-1,\ldots,1}{n,n-1,\ldots,0}$.
Therefore, by the same argument in the symplectic case,
we obtain $\Theta\text{\rm -rk}(\rho)=\Theta\text{\rm -rk}(\rho\cdot\sgn)=2n$.

\item Suppose that $\bfG=\rmU_n$ and $n\geq 1$.
\begin{enumerate}
\item Suppose that $n$ is even.
The Steinberg character of $G$ is the unipotent character associated to the symbol
$\Lambda=\binom{\frac{n}{2}-1,\frac{n}{2}-2,\ldots,0}{\frac{n}{2},\frac{n}{2}-1,\ldots,1}$
and $\Upsilon(\Lambda)=\sqbinom{-}{1,1,\ldots,1}$ ($\frac{n}{2}$ copies of ``$1$'').
Then
\[
\Theta\text{\rm -rk}(\St_{\rmU_n})
={\rm rk}_\rmU\textstyle\binom{\frac{n}{2}-1,\frac{n}{2}-2,\ldots,1}{\frac{n}{2}-1,\frac{n}{2}-2,\ldots,0}
=2\left(\frac{n}{2}-1\right)+\frac{1\cdot 2}{2}=n-1.
\]

\item Suppose that $n$ is odd.
The Steinberg character of $G$ is the unipotent character associated to the symbol
$\Lambda=\binom{\frac{n-1}{2},\frac{n-3}{2},\ldots,1}{\frac{n-1}{2},\frac{n-3}{2},\ldots,0}$
and $\Upsilon(\Lambda)=\sqbinom{1,1,\ldots,1}{-}$ ($\frac{n-1}{2}$ copies of ``$1$'').
Then
\[
\Theta\text{\rm -rk}(\St_{\rmU_n})
={\rm rk}_\rmU\textstyle\binom{\frac{n-1}{2},\frac{n-3}{2},\ldots,0}{\frac{n-3}{2},\frac{n-5}{2},\ldots,1}
=2\left(\frac{n-1}{2}\right)+\frac{0\cdot 1}{2}=n-1.
\]
\end{enumerate}
\end{enumerate}
\end{exam}

\begin{lem}\label{0309}
Let $\bfG$ be $\Sp_{2n},\rmO^+_{2n}$, or $\rmO_{2n+1}$, and let $\rho\in\cale(G)_1$.
If\/ $\Theta\text{\rm -rk}(\rho)=2n$, then
\[
\begin{cases}
\text{$\rho$ is a Steinberg character}, & \text{if\/ $\bfG=\Sp_{2n},\rmO_{2n+1}$};\\
n=0, & \text{if\/ $\bfG=\rmO^+_{2n}$}.
\end{cases}
\]
Moreover, $\rmO^-_{2n}$ does not have any unipotent character of\/ $\Theta$-rank $2n$.
\end{lem}
\begin{proof}
Let $\rho=\rho_\Lambda\in\cale(G)_1$
($\rho=\rho_\Lambda$ or $\rho=\rho_\Lambda\cdot\sgn$ if $\bfG=\rmO_{2n+1}$)
such that $\Theta\text{\rm -rk}(\rho)=2n$,
and write
\[
\Lambda=\binom{a_1,a_2,\ldots,a_{m_1}}{b_1,b_2,\ldots,b_{m_2}}\in\cals_\bfG,\qquad
\Upsilon(\Lambda)=\sqbinom{\mu_1,\mu_2,\ldots,\mu_{m_1}}{\nu_1,\nu_2,\ldots,\nu_{m_2}}.
\]
\begin{enumerate}
\item Suppose that $\bfG=\Sp_{2n}$.
By Proposition~\ref{0306}, we need both ${\rm rk}\binom{b_1,b_2,\ldots,b_{m_2}}{a_2,a_3,\ldots,a_{m_1}}=n$
and ${\rm rk}\binom{b_2,b_3,\ldots,b_{m_2}}{a_1,a_2,\ldots,a_{m_1}}=n$,
and then $\binom{b_1,b_2,\ldots,b_{m_2}}{a_2,a_3,\ldots,a_{m_1}}\in\cals_{\rmO^+_{2n}}$
and $\binom{b_2,b_3,\ldots,b_{m_2}}{a_1,a_2,\ldots,a_{m_1}}\in\cals_{\rmO^-_{2n}}$.
So we need $\mu_1=0$ and $\nu_1=1$, thus
$\Upsilon(\Lambda)=\sqbinom{-}{1,1,\ldots,1}$ ($n$ copies of ``$1$'').
That is, $\rho_\Lambda$ is the Steinberg character of $\Sp_{2n}(q)$.

\item Suppose that $\bfG=\rmO_{2n+1}$.
By Corollary~\ref{0310}, we know that $\Theta\text{\rm -rk}(\rho_\Lambda)=2n$ if and only if
$\Theta\text{\rm -rk}(\rho_{\Lambda^\rmt})=2n$.
Now by (1), $\rho_{\Lambda^\rmt}$ is the Steinberg character of $\Sp_{2n}(q)$,
hence $\rho_\Lambda$ is a Steinberg character of $\rmO_{2n+1}(q)$.

\item Suppose that $\bfG=\rmO^+_{2n}$.
By Proposition~\ref{0306}, we need both ${\rm rk}\binom{b_1,b_2,\ldots,b_{m_2}}{a_2,a_3,\ldots,a_{m_1}}=n$
and ${\rm rk}\binom{b_2,b_3,\ldots,b_{m_2}}{a_1,a_2,\ldots,a_{m_1}}=n$,
and then $\binom{b_1,b_2,\ldots,b_{m_2}}{a_2,a_3,\ldots,a_{m_1}}\in\cals_{\Sp_{2n}}$
and $\binom{a_1,a_2,\ldots,a_{m_1}}{b_2,b_3,\ldots,b_{m_2}}\in\cals_{\Sp_{2n}}$.
So we need $\mu_1=0$ and $\nu_1=0$, thus
$\Upsilon(\Lambda)=\sqbinom{-}{-}$, i.e., $\rho_\Lambda$ is a cuspidal unipotent character.
But from Example~\ref{0317}, we must have $n=0$.

\item Suppose that $\bfG=\rmO^-_{2n}$ (where $n\geq 1$).
Now $|\Upsilon(\Lambda)|=n-1$ and the minimum value of
${\rm rk}(\Lambda^\rmt\smallsetminus\binom{-}{a_1})$ and ${\rm rk}(\Lambda^\rmt\smallsetminus\binom{b_1}{-})$
is less than or equal to $n-1$.
This means that the maximum value of $\Theta\text{\rm -rk}(\rho)$ for $\rho\in\cale(\rmO^-_{2n}(q))_1$
is at most $2(n-1)$.
\end{enumerate}
\end{proof}

\begin{lem}\label{0313}
Let $\rho$ be an irreducible pseudo-unipotent character of\/ $\Sp_{2n}(q)$
and write $\Xi_s(\rho)=\rho^{(0)}\otimes\rho_{\Lambda^{(-)}}\otimes\rho_{\Lambda^{(+)}}$.
Then
\[
\Theta\text{\rm -rk}(\rho)=\Theta\text{\rm -rk}(\rho_{\Lambda^{(-)}})+1,
\]
in particular, $\Theta\text{\rm -rk}(\rho)$ is odd.
\end{lem}
\begin{proof}
Let $\rho$ be an irreducible pseudo-unipotent character of $G=\Sp_{2n}(q)$.
Then by Lemma~\ref{0209}, to determine $\Theta\text{\rm -rk}(\rho)$,
we only need to consider the first occurrence of $\rho$ for the dual pair $(\bfG',\bfG)=(\rmO_{2n'+1},\Sp_{2n})$.
Now we have $\rho\in\cale(G)_s$ such that
$\bfG^{(0)}=\rmU_0$, $\bfG^{(-)}=\rmO^{\epsilon'}_{2n}$ for some $\epsilon'$, and $\bfG^{(+)}=\Sp_0$.
Consider the diagram
\[
\begin{CD}
\rho' @>>> \rho \\
@VVV @VVV \\
{\bf 1}\otimes\rho_{\Lambda'^{(+)}}\otimes{\bf 1} @>>> {\bf 1}\otimes\rho_{\Lambda^{(-)}}\otimes{\bf 1}
\end{CD}
\]
Then $\rho'\in\cale(G')_{s'}$ such that $\bfG'^{(0)}=\rmU_0$, $\bfG'^{(-)}=\Sp_{2n'}$ and $\bfG^{(+)}=\Sp_0$.
Now $\rho$ of $\Sp_{2n}(q)$ occurs in the correspondence for the dual pair $(\rmO_{2n'+1},\Sp_{2n})$
if and only if $\rho_{\Lambda^{(-)}}$ of $\rmO^{\epsilon'}_{2n}(q)$ occurs in the correspondence for the dual
pair $(\Sp_{2n'},\rmO^{\epsilon'}_{2n})$.
Thus the lemma is proved.
\end{proof}

\begin{lem}\label{0312}
Let $\rho\in\cale(\rmU_n(q))_1$.
If\/ $\Theta\text{\rm -rk}(\rho)=n$,
then $n=0$.
\end{lem}
\begin{proof}
Suppose that $\rho=\rho_\Lambda\in\cale(G)_1$.
By the same argument in the proof of Lemma~\ref{0309},
we must have $\Upsilon(\Lambda)=\sqbinom{0}{0}$, i.e., $\rho_\Lambda$ is cuspidal.
This means that ${\rm rk}_\rmU(\Lambda)=\frac{1}{2}d(d+1)$ and $\Theta\text{\rm -rk}(\Lambda)=\frac{1}{2}d(d-1)$
for some non-negative integer $d$.
The equality $n=\frac{1}{2}d(d+1)=\frac{1}{2}d(d-1)$ implies that $n=d=0$.
\end{proof}

\subsection{$\Theta$-ranks for symplectic groups or orthogonal groups}
Now we can describe the $\Theta$-rank of a general irreducible character of a symplectic group
or an orthogonal group in terms of its data under the modified Lusztig correspondence.

\begin{prop}\label{0307}
Let $\bfG$ be $\Sp_{2n}$, $\rmO^\epsilon_{2n}$, or $\rmO_{2n+1}$, and let $\rho\in\cale(G)$.
Suppose that $\rho\in\cale(G)_s$ for some $s$ and write
$\Xi_s(\rho)=\rho^{(0)}\otimes\rho_{\Lambda^{(-)}}\otimes\rho_{\Lambda^{(+)}}$.
Then $\Theta\text{\rm -rk}(\rho)$ is equal to
\[
\begin{cases}
\min\{2n-2n^{(+)}+\Theta\text{\rm -rk}(\Lambda^{(+)}),2n-2n^{(-)}+\Theta\text{\rm -rk}(\Lambda^{(-)})+1\},
& \text{if\/ $\bfG=\Sp_{2n}$};\\
\min\{2n-2n^{(+)}+\Theta\text{\rm -rk}(\Lambda^{(+)}),2n-2n^{(-)}+\Theta\text{\rm -rk}(\Lambda^{(-)})\},
& \text{if\/ $\bfG=\rmO^\epsilon_{2n},\rmO_{2n+1}$}
\end{cases}
\]
where $n^{(-)},n^{(+)}$ are given in Subsection~\ref{0215},
and $\Theta\text{\rm -rk}(\Lambda^{(+)}), \Theta\text{\rm -rk}(\Lambda^{(-)})$ are given in (\ref{0201}).
\end{prop}
\begin{proof}
First suppose that $\rho'\otimes\rho$ occurs in the correspondence for the dual pair
$(\bfG',\bfG)$ for some $\rho'\in\cale(G')_{s'}$ and write $\Xi_{s'}(\rho')=\rho'^{(0)}\otimes\rho_{\Lambda'^{(-)}}\otimes\rho_{\Lambda'^{(+)}}$.
\begin{enumerate}
\item Suppose that $\bfG=\Sp_{2n}$.
First we consider the dual pair $(\bfG',\bfG)=(\rmO^\epsilon_{2n'},\Sp_{2n})$ for some $n'$.
Then we have a commutative diagram
\[
\begin{CD}
\rho' @> \Theta_{\bfG',\bfG} >> \rho \\
@V \Xi_{s'} VV @VV \Xi_s V \\
\rho'^{(0)}\otimes\rho_{\Lambda'^{(-)}}\otimes\rho_{\Lambda'^{(+)}} @> \id\otimes\id\otimes\Theta_{\bfG'^{(+)},\bfG^{(+)}} >>  \rho^{(0)}\otimes\rho_{\Lambda^{(-)}}\otimes\rho_{\Lambda^{(+)}}
\end{CD}
\]
Now we have
$\bfG'^{(-)}=\rmO^{\epsilon'\epsilon}_{2n'^{(-)}}$ and
$\bfG'^{(+)}=\rmO^{\epsilon'}_{2n'^{(+)}}$;
$\bfG^{(-)}=\rmO^{\epsilon'\epsilon}_{2n^{(-)}}$ and
$\bfG^{(+)}=\Sp_{2n^{(+)}}$ for some $\epsilon'$ depending on $s$.
Because now $G'^{(0)}=G^{(0)}$ and $G'^{(-)}=G^{(-)}$,
we have $2n'-2n'^{(+)}=2n-2n^{(+)}$, i.e.,
$2n'=2n-2n^{(+)}+2n'^{(+)}$.
Now $n$ and $n^{(+)}$ are fixed and the minimal value of $2n'^{(+)}$ (for both possible $\epsilon=+$ or $-$)
is equal to $\Theta\text{\rm -rk}(\Lambda^{(+)})$.

Next we consider the dual pair $(\bfG',\bfG)=(\rmO_{2n'+1},\Sp_{2n})$.
Then we have a commutative diagram
\[
\begin{CD}
\rho' @> \Theta_{\bfG',\bfG} >> \rho \\
@V \iota\circ\Xi_{s'} VV @VV \Xi_s V \\
\rho'^{(0)}\otimes\rho_{\Lambda'^{(+)}}\otimes\rho_{\Lambda'^{(-)}} @> \id\otimes\Theta_{\bfG'^{(+)},\bfG^{(-)}}\otimes\id
>> \rho^{(0)}\otimes\rho_{\Lambda^{(-)}}\otimes\rho_{\Lambda^{(+)}}
\end{CD}
\]
Now we have
$\bfG'^{(-)}=\Sp_{2n'^{(-)}}$ and
$\bfG'^{(+)}=\Sp_{2n'^{(+)}}$;
$\bfG^{(-)}=\rmO^{\epsilon'}_{2n^{(-)}}$ and
$\bfG^{(+)}=\Sp_{2n^{(+)}}$ for some $\epsilon'$.
Because now $G'^{(0)}=G^{(0)}$ and $G'^{(-)}=G^{(+)}$,
we have $2n'-2n'^{(+)}=2n-2n^{(-)}$, i.e.,
$2n'+1=2n-2n^{(-)}+2n'^{(+)}+1$.
Now $n$ and $n^{(-)}$ are fixed and the minimal value of $2n'^{(+)}$
is equal to $\Theta\text{\rm -rk}(\Lambda^{(-)})$.

From the definition of $\Theta$-rank and Proposition~\ref{0306}, we see that
\begin{align*}
\Theta\text{\rm -rk}(\rho)
&=\min\{2n-2n^{(+)}+\Theta\text{\rm -rk}(\rho_{\Lambda^{(+)}}),
2n-2n^{(-)}+\Theta\text{\rm -rk}(\rho_{\Lambda^{(-)}})+1\} \\
&=\min\{2n-2n^{(+)}+\Theta\text{\rm -rk}(\Lambda^{(+)}),
2n-2n^{(-)}+\Theta\text{\rm -rk}(\Lambda^{(-)})+1\}.
\end{align*}

\item Suppose that $\bfG=\rmO^\epsilon_{2n}$.
Now we consider the dual pair $(\Sp_{2n'},\rmO^\epsilon_{2n})$ and we have a commutative diagram
\[
\begin{CD}
\rho' @> \Theta >> \rho \\
@V \Xi_{s'} VV @VV \Xi_s V \\
\rho'^{(0)}\otimes\rho_{\Lambda'^{(-)}}\otimes\rho_{\Lambda'^{(+)}} @> \id\otimes\id\otimes\Theta >>  \rho^{(0)}\otimes\rho_{\Lambda^{(-)}}\otimes\rho_{\Lambda^{(+)}}
\end{CD}
\]
Now we have
$\bfG'^{(-)}(s')=\rmO^{\epsilon'\epsilon}_{2n'^{(-)}}$ and
$\bfG'^{(+)}(s')=\Sp_{2n'^{(+)}}$;
$\bfG^{(-)}(s)=\rmO^{\epsilon'\epsilon}_{2n^{(-)}}$ and
$\bfG^{(+)}(s)=\rmO^{\epsilon'}_{2n^{(+)}}$ for some $\epsilon'$.

Now we need consider the occurrences for four irreducible characters
$\rho$, $\rho\cdot\sgn$, $\rho\chi_\bfG$ and $\rho\chi_\bfG\cdot\sgn$.
Note that by Lemma~\ref{0207}, we have
$\Lambda^{(\varepsilon)}(\rho\cdot\sgn)=(\Lambda^{(\varepsilon)}(\rho))^\rmt$
and $\Lambda^{(\varepsilon)}(\rho\chi_\bfG)=\Lambda^{(-\varepsilon)}(\rho)$ for $\varepsilon=+$ or $-$.
By the similar argument in (1), we conclude that
\[
\Theta\text{\rm -rk}(\rho)
=\min\{2n-2n^{(+)}+\Theta\text{\rm -rk}(\Lambda^{(+)}),
2n-2n^{(-)}+\Theta\text{\rm -rk}(\Lambda^{(-)})\}.
\]

\item Suppose that $\bfG=\rmO_{2n+1}$.
Now consider the dual pair $(\Sp_{2n'},\rmO_{2n+1})$ and we have a commutative diagram
\[
\begin{CD}
\rho' @> \Theta >> \rho \\
@V \Xi_{s'} VV @VV \Xi_s V \\
\rho'^{(0)}\otimes\rho_{\Lambda'^{(-)}}\otimes\rho_{\Lambda'^{(+)}} @> \id\otimes\Theta\otimes\id >>  \rho^{(0)}\otimes\rho_{\Lambda^{(-)}}\otimes\rho_{\Lambda^{(+)}}
\end{CD}
\]
Now we have
$\bfG'^{(-)}(s')=\rmO^{\epsilon'\epsilon}_{2n'^{(-)}}$ and
$\bfG'^{(+)}(s')=\Sp_{2n'^{(+)}}$;
$\bfG^{(-)}(s)=\Sp_{2n^{(-)}}$ and
$\bfG^{(+)}(s)=\Sp_{2n^{(+)}}$ for some $\epsilon'$.
Then by the similar argument in (2), we conclude that
\[
\Theta\text{\rm -rk}(\rho)
=\min\{2n-2n^{(+)}+\Theta\text{\rm -rk}(\Lambda^{(+)}),
2n-2n^{(-)}+\Theta\text{\rm -rk}(\Lambda^{(-)})\}.
\]
\end{enumerate}
\end{proof}

\begin{exam}
There are two irreducible characters $\rho_1,\rho_2$ of degree $\frac{1}{2}(q^n+1)$ of $\Sp_{2n}(q)$
for $n\geq 1$.
They are in the Lusztig series $\cale(G)_s$ such that
$\bfG^{(0)}=\rmU_0$, $\bfG^{(-)}=\rmO^+_{2n}$, $\bfG^{(+)}=\Sp_0$, i.e., $n^{(-)}=n$ and $n^{(+)}=0$.
Moreover, we know that $\Lambda^{(-)}=\binom{n}{0}$ or $\binom{0}{n}$,
and $\Lambda^{(+)}=\binom{0}{-}$.
Therefore
\begin{align*}
\Theta\text{\rm -rk}(\rho_1)=\Theta\text{\rm -rk}(\rho_2)
&=\min\{2n-0+\Theta\text{\rm -rk}(\textstyle\binom{0}{-}),2n-2n+\Theta\text{\rm -rk}(\textstyle\binom{n}{0})+1\} \\
&=\min\{2n,1\}=1.
\end{align*}
In fact, both $\rho_1,\rho_2$ are pseudo-unipotent,
and each of $(\rho_1)_{\Lambda^{(-)}},(\rho_2)_{\Lambda^{(-)}}$ is the trivial character or the $\sgn$
character of $\rmO^+_{2n}(q)$.
So the above result is also a direct consequence of Lemma~\ref{0313}.
\end{exam}

\begin{cor}
Let $\bfG=\Sp_{2n},\rmO_{2n+1}$, or $\rmO^\epsilon_{2n}$,
and let $\rho\in\cale(G)_s$ for some $s$.
Then $\Theta\text{\rm -rk}(\rho)=2n$ if and only if
\[
\begin{cases}
\text{$n^{(-)}=0$ and $\rho_{\Lambda^{(+)}}$ is a Steinberg character}, & \text{if\/ $\bfG=\Sp_{2n}$};\\
\text{both $\rho_{\Lambda^{(-)}},\rho_{\Lambda^{(+)}}$ are Steinberg characters}, & \text{if\/ $\bfG=\rmO_{2n+1}$};\\
\text{both $n^{(-)}=n^{(+)}=0$}, & \text{if\/ $\bfG=\rmO^\epsilon_{2n}$}.
\end{cases}
\]
\end{cor}
\begin{proof}
Suppose that $\bfG=\Sp_{2n}$.
By Proposition~\ref{0307},
we know that $\Theta\text{\rm -rk}(\rho)=2n$ if and only if $\Theta\text{\rm -rk}(\rho_{\Lambda^{(-)}})=2n^{(-)}$
and $\Theta\text{\rm -rk}(\rho_{\Lambda^{(+)}})=2n^{(+)}$.
Now $\bfG^{(-)}=\rmO_{2n^{(-)}}$ and $G^{(+)}=\Sp_{2n^{(+)}}$,
so we need $n^{(-)}=0$ and $\rho^{(+)}$ to be the Steinberg character by Lemma~\ref{0309}.
The proofs for two other cases are similar.
\end{proof}

We provide some necessary condition in Remark~\ref{0304} on the possible values of
$\Theta\text{\rm -rk}(\rho)$ for $\rho\in\cale(G)$.
Now we also have the sufficient condition:

\begin{prop}
Let $\bfG=\rmO^\epsilon_{2n}$, $\Sp_{2n}$ or $\rmO_{2n+1}$ over $\bff_q$.
\begin{enumerate}
\item[(i)] Suppose that $\bfG=\rmO^\epsilon_{2n}$.
\begin{enumerate}
\item For each even integer $k$ such that $0\leq k<2n$,
there exists an irreducible character $\rho\in\cale(\rmO^\epsilon_{2n}(q))$ such that $\Theta\text{\rm -rk}(\rho)=k$.

\item If $2n\geq 4$, or $q>3$, or $\epsilon=-$,
then there exists an irreducible character $\rho\in\cale(\rmO^\epsilon_{2n}(q))$ such that
$\Theta\text{\rm -rk}(\rho)=2n$.

\item There is no irreducible character of\/ $\rmO^+_2(3)$ of\/ $\Theta$-rank $2$.
\end{enumerate}

\item[(ii)] For each integer $k$ such that $0\leq k\leq 2n$,
there exists an irreducible character $\rho\in\cale(\Sp_{2n}(q))$ such that $\Theta\text{\rm -rk}(\rho)=k$.

\item[(iii)] For each even integer $k$ such that $0\leq k\leq 2n$,
there exists an irreducible character $\rho\in\cale(\rmO_{2n+1}(q))$ such that $\Theta\text{\rm -rk}(\rho)=k$.
\end{enumerate}
\end{prop}
\begin{proof}
\begin{enumerate}
\item Suppose that $\bfG=\rmO^-_{2n}$.
The assertion is trivial when $k=0$,
so we will assume that $k$ is an even integer such that $2\leq k\leq 2n$.
\begin{enumerate}
\item Suppose that $2\leq k<2n$.
Let $\Lambda_k=\binom{n-1,\frac{k}{2},\frac{k}{2}-1,\ldots,1}{\frac{k}{2}-2,\frac{k}{2}-3,\ldots,0}\in\cals_{\rmO^-_{2n}}$.
Then we have $\Upsilon(\Lambda_k)=\sqbinom{-}{n-\frac{k}{2}-1,1,\ldots,1}\in\calp_2(n-1)$.
Then $\binom{\frac{k}{2},\frac{k}{2}-1,\ldots,1}{\frac{k}{2}-2,\frac{k}{2}-3,\ldots,0}\in\cals_{\Sp_k}$.
Moreover, $\binom{-}{n,2,0}\in\cals_{\Sp_{2n+2}}$ if $k=2$;  $\binom{\frac{k}{2}-3,\frac{k}{2}-4,\ldots,0}{n-1,\frac{k}{2},\frac{k}{2}-1,\ldots,1}\in\cals_{\Sp_{2n+2}}$
if $k\geq 4$.
Then by Proposition~\ref{0306}, 
$\Theta\text{\rm -rk}(\rho_{\Lambda_k})=\min\{k,2n+2\}=k$.

\item It is known that there exists a semisimple element $s\in G^*=\rmO^-_{2n}(q)$ such that
$C_{G^*}(s)\simeq\rmU_1(q^n)$.
That is, $n^{(-)}=n^{(+)}=0$, and then any $\rho\in\cale(G)_s$ has $\Theta\text{\rm -rk}(\rho)=2n$ by Proposition~\ref{0307}.
\end{enumerate}

\item Suppose that $\bfG=\rmO_{2n}^+$ and $k$ is an even integer such that $0\leq k\leq 2n$.

\begin{enumerate}
\item Suppose that $k<2n$.
Let $\Lambda_k=\binom{n,\frac{k}{2},\frac{k}{2}-1,\ldots,1}{\frac{k}{2},\frac{k}{2}-1,\ldots,0}\in\cals_{\rmO^+_{2n}}$.
Then we have $\Upsilon(\Lambda_k)=\sqbinom{n-\frac{k}{2},1,\ldots,1}{-}\in\calp_2(n)$,
$\binom{\frac{k}{2},\frac{k}{2}-1,\ldots,0}{\frac{k}{2},\frac{k}{2}-1,\ldots,1}\in\cals_{\Sp_{k}}$,
and $\binom{n,\frac{k}{2},\frac{k}{2}-1,\ldots,1}{\frac{k}{2}-1,\frac{k}{2}-2,\ldots,0}\in\cals_{\Sp_{2n}}$.
Hence by Proposition~\ref{0306}, $\Theta\text{\rm -rk}(\rho_{\Lambda_k})=\min\{k,2n\}=k$.

\item Suppose that $q>3$ or $n\geq 2$.
It is known that there exists a semisimple element $s\in G^*=\rmO^+_{2n}(q)$ such that
$C_{G^*}(s)\simeq\GL_1(q^n)$.
That is, $n^{(-)}=n^{(+)}=0$, and then any $\rho\in\cale(G)_s$ has $\Theta\text{\rm -rk}(\rho)=2n$ by Proposition~\ref{0307}.

\item Now $\rmO^+_2(3)$ is a finite group of $4$ element.
Each irreducible character of $\rmO^+_2(3)$ is a linear character and hence has $\Theta$-rank $0$.
\end{enumerate}

\item Suppose that $\bfG=\Sp_{2n}$ and $k$ is an integer such that $0\leq k\leq 2n$.
\begin{enumerate}
\item Suppose that $k$ is even and $0\leq k<2n$.
Let $\Lambda_k=\binom{n,\frac{k}{2},\frac{k}{2}-1,\ldots,1}{\frac{k}{2}-1,\frac{k}{2}-2,\ldots,0}\in\cals_{\Sp_{2n}}$ 
and then $\Upsilon(\Lambda_k)=\sqbinom{n-\frac{k}{2},1,\ldots,1}{-}\in\calp_2(n)$.
Now $\binom{\frac{k}{2}-1,\frac{k}{2}-2,\ldots,0}{\frac{k}{2},\frac{k}{2}-1,\ldots,1}\in\cals_{\rmO^+_{2k}}$, $\binom{\frac{k}{2}-2,\frac{k}{2}-3,\ldots,0}{n,\frac{k}{2},\frac{k}{2}-1,\ldots,1}\in\cals_{\rmO^-_{2n+2}}$.
Thus by Proposition~\ref{0306}, 
$\Theta\text{\rm -rk}(\rho_{\Lambda_k})=\min\{k,2n+2\}=k$.

\item Suppose that $k=2n$.
Let $\Lambda_k=\binom{n,n-1,\ldots,0}{n,n-1,\ldots,1}\in\cals_{\Sp_{2n}}$.
Then $\Upsilon(\Lambda_k)=\sqbinom{-}{1,1,\ldots,1}\in\calp_2(n)$.
Now $\binom{n,n-1,\ldots,1}{n-1,n-2,\ldots,0}\in\cals_{\rmO^+_{2n}}$ and
$\binom{n-1,n-2,\ldots,1}{n,n-1,\ldots,0}\in\cals_{\rmO^-_{2n}}$.
Thus $\Theta\text{\rm -rk}(\rho_{\Lambda_k})=2n$.

\item Suppose that $k$ is odd.
Then $k-1$ is even and $0\leq k-1<2n$.
Let $s$ be an semisimple element in $G^*$ such that $\bfG^{(0)}=\rmU_0$, $\bfG^{(-)}=\rmO^+_{2n}$ and $\bfG^{(+)}=\Sp_0$.
Let $\rho\in\cale(G)_s$ such that $\Xi_s(\rho)={\bf 1}\otimes\rho_{\Lambda^{(-)}}\otimes{\bf 1}$ where
$\Lambda^{(-)}\in\cals_{\rmO^+_{2n}}$ is the symbol $\Lambda_{k-1}$ given in (2.a).
Now $\rho$ is pseudo-unipotent (\cf.~Subsection~\ref{0204}), 
and by Lemma~\ref{0313}, $\Theta\text{\rm -rk}(\rho)=\Theta\text{\rm -rk}(\rho_{\Lambda_{k-1}})+1=k-1+1=k$.
\end{enumerate}

\item Suppose that $\bfG=\rmO_{2n+1}$ and $k$ is an even integer such that $0\leq k\leq 2n$.
Let $\rho=\rho_\Lambda$ where $\Lambda^\rmt$ is equal to the symbol $\Lambda_k\in\cals_{\Sp_{2n}}$
given in (3.a) and (3.b).
Then by Corollary~\ref{0310}, we have $\Theta\text{\rm -rk}(\rho)=\Theta\text{\rm -rk}(\rho_{\Lambda^\rmt})=k$.
\end{enumerate}
\end{proof}

\subsection{$\Theta$-ranks for unitary groups}\label{0316}
In this subsection, we consider $\bfG=\rmU_n$.
For $s\in G^*$, we can write
\[
C_{G^*}(s)=\prod_{i=1}^r\prod_{j=1}^{t_i}\GL_{n_{ij}}^{(-1)^i}(q^i)
\]
for some non-negative integers $r,t_i,n_{ij}$ such that $\sum_{i=1}^r\sum_{j=1}^{t_i}in_{ij}=n$ where
$\GL^{+1}_{n_{ij}}=\GL^{+}_{n_{ij}}=\GL_{n_{ij}}$
and $\GL^{-1}_{n_{ij}}=\GL^{-}_{n_{ij}}=\rmU_{n_{ij}}$.
For $i=1$, we can let $t_1=q+1$ and each $j=1,\ldots,t_1$ corresponds an
eigenvalue $\lambda_j\in\overline\bff_q$ of $s$ such that $\lambda_j^{q+1}=1$.
Then for the Lusztig correspondence $\grL_s\colon\cale(G)_1\rightarrow\cale(C_{G^*}(s))_1$
we can write
\begin{equation}\label{0308}
\grL_s(\rho)=\bigotimes_{i=1}^r\bigotimes_{j=1}^{t_i}\rho^{(ij)}
\end{equation}
for some $\rho^{(ij)}\in\cale(\GL_{n_{ij}}^{(-1)^i}(q^i))_1$.
Note that $\rho^{(1,j)}$ is a unipotent character of $\rmU_{n_{1,j}}(q)$.

\begin{prop}\label{0311}
Let $\rho\in\cale(\rmU_n(q))_s$ for some $s$, and keep the above notations.
Then
\[
\Theta\text{\rm -rk}(\rho)=\min\{n-n_{1,j}+\Theta\text{\rm -rk}(\rho^{(1,j)})\mid j=1,\ldots,q+1\}.
\]
\end{prop}
\begin{proof}
Write $\grL_s(\rho)$ as in (\ref{0308}).
Suppose that $\bfG^{(+)}(s)=\rmU_{n_{1,1}}$.
Then the minimal value of $n'$ such that $\rho$ occurs in the $\Theta$-correspondence for the pair
$(\rmU_{n'},\rmU_n)$ is equal to $n-n_{1,1}+\Theta\text{\rm -rk}(\rho^{(1,1)})$.
Suppose that $\chi$ is a linear character of $\rmU_n(q)$.
Then $\rho\chi\in\cale(G)_{s'}$ where $s'=sa$ for some $a\in\bff_{q^2}$ such that $a^{q+1}=1$.
Then we can write
\[
C_{G^*}(s')=\prod_{i=1}^r\prod_{j=1}^{t_i}\GL_{n'_{ij}}^{(-1)^i}(q^i)\quad\text{and}\quad
\grL_{s'}(\rho\chi)=\bigotimes_{i=1}^r\bigotimes_{j=1}^{t_i}\rho'^{(ij)}
\]
where $\{\rmU_{n'_{1,1}},\ldots\rmU_{n'_{1,t_1}}\}$ is a permutation of $\{\rmU_{n_{1,1}},\ldots\rmU_{n_{1,t_1}}\}$ and
$\{\rho'^{(1,1)},\ldots,\rho'^{(1,t_1)}\}$ are a permutation of $\{\rho^{(1,1)},\ldots,\rho^{(1,t_1)}\}$.
Moreover, for any $1\leq j\leq t_1$ there exists a linear character $\chi$ of $\rmU_n(q)$ such that
$\rho'^{(1,1)}=\rho^{(1,j)}$.
Hence the proposition follows from the definition of the $\Theta$-rank.
\end{proof}

\begin{rem}
The same result of this proposition can be found in \cite{GLT01} theorem 3.9.
\end{rem}

\begin{cor}\label{0314}
Let $\rho\in\cale(\rmU_n(q))$.
Then $\Theta\text{\rm -rk}(\rho)=n$ if and only if $n_{1,j}=0$ for each $j=1,\ldots,q+1$.
\end{cor}
\begin{proof}
By Proposition~\ref{0311}, $\Theta\text{\rm -rk}(\rho)=n$ implies that
$\Theta\text{\rm -rk}(\rho^{(1,j)})=n_{1,j}$ for each $j=1,\ldots,q+1$.
Then by Lemma~\ref{0312}, we conclude that $n_{1,j}=0$ for each $j=1,\ldots,q+1$.
\end{proof}

\begin{prop}
For each integer $k$ such that $k=0$ if $n=0,1$; or $0\leq k\leq n$ if $n\geq 2$,
there exists an irreducible character $\rho\in\cale(\rmU_n(q))$ such that $\Theta\text{\rm -rk}(\rho)=k$.
In particular, there is no $\rho\in\cale(\rmU_1(q))$ of $\Theta$-rank $1$.
\end{prop}
\begin{proof}
The case for $n=0$ is trivial.
For $n=1$, every irreducible character of $\rmU_1(q)$ is a linear character and hence 
of $\Theta$-rank $0$.
Now we assume $n\geq 2$.
\begin{enumerate}
\item Suppose that $0\leq k<n$.

\begin{enumerate}
\item Suppose that $n+k$ is even.
Let
\[
\Lambda_k=\begin{cases}
\binom{\frac{n}{2},\frac{k}{2},\frac{k}{2}-1,\ldots,1}{\frac{k}{2},\frac{k}{2}-1,\ldots,0};\\
\binom{\frac{n-1}{2},\frac{k-1}{2},\frac{k-3}{2},\ldots,1}{\frac{k+1}{2},\frac{k-1}{2},\ldots,0},
\end{cases}\quad
\Upsilon(\Lambda_k)=\begin{cases}
\sqbinom{\frac{n-k}{2},1,1,\ldots,1}{-}, & \text{if $n$ is even;}\\
\sqbinom{\frac{n-k}{2},1,1,\ldots,1}{-}, & \text{if $n$ is odd.}
\end{cases}
\]
\begin{enumerate}
\item Suppose that $n$ is even.
Now $\binom{\frac{k}{2},\frac{k}{2}-1,\ldots,0}{\frac{k}{2},\frac{k}{2}-1,\ldots,1}\in\cals_{\rmU_k}$ and
$\binom{\frac{k}{2}-1,\frac{k}{2}-2,\ldots,0}{\frac{n}{2},\frac{k}{2},\frac{k}{2}-1,\ldots,1}\in\cals_{\rmU_{n+1}}$.
Hence $\Theta\text{\rm -rk}(\rho_{\Lambda_k})=\min\{k,n+1\}=k$.

\item Suppose that $n$ is odd.
Now $\binom{\frac{k+1}{2},\frac{k-1}{2},\ldots,0}{\frac{k-1}{2},\frac{k-3}{2},\ldots,1}\in\cals_{\rmU_k}$ and
$\binom{\frac{k-1}{2},\frac{k-3}{2},\ldots,0}{\frac{n-1}{2},\frac{k-1}{2},\frac{k-3}{2},\ldots,1}\in\cals_{\rmU_{n-1}}$.
Hence $\Theta\text{\rm -rk}(\rho_{\Lambda_k})=\min\{k,n-1\}=k$.
\end{enumerate}

\item Suppose that $n+k$ is odd.
Let
\[
\Lambda_k=\begin{cases}
\binom{\frac{k-1}{2},\frac{k-3}{2},\ldots,0}{\frac{n}{2},\frac{k-1}{2},\frac{k-3}{2},\ldots,1};\\
\binom{\frac{k}{2}-1,\frac{k}{2}-2,\ldots,0}{\frac{n-1}{2},\frac{k}{2},\frac{k}{2}-1,\ldots,1},
\end{cases}\quad
\Upsilon(\Lambda_k)=\begin{cases}
\sqbinom{-}{\frac{n-k+1}{2},1,1,\ldots,1}, & \text{if $n$ is even;}\\
\sqbinom{-}{\frac{n-k-1}{2},1,1,\ldots,1}, & \text{if $n$ is odd.}
\end{cases}
\]
\begin{enumerate}
\item Suppose that $n$ is even.
Now $\binom{\frac{k-1}{2},\frac{k-3}{2},\ldots,1}{\frac{k-1}{2},\frac{k-3}{2},\ldots,0}\in\cals_{\rmU_k}$ and
$\binom{\frac{n}{2},\frac{k-1}{2},\frac{k-3}{2},\ldots,1}{\frac{k-1}{2},\frac{k-3}{2},\ldots,0}\in\cals_{\rmU_n}$.
Hence $\Theta\text{\rm -rk}(\rho_{\Lambda_k})=\min\{k,n\}=k$.

\item Suppose that $n$ is odd.
Now $\binom{\frac{k}{2},\frac{k}{2}-1,\ldots,1}{\frac{k}{2}-1,\frac{k}{2}-2,\ldots,0}\in\cals_{\rmU_k}$ and
$\binom{\frac{n-1}{2},\frac{k}{2},\frac{k}{2}-1,\ldots,1}{\frac{k}{2}-2,\frac{k}{2}-3,\ldots,0}\in\cals_{\rmU_{n+2}}$.
Hence $\Theta\text{\rm -rk}(\rho_{\Lambda_k})=\min\{k,n+2\}=k$.
\end{enumerate}
\end{enumerate}

\item Suppose that $k=n$.
\begin{enumerate}
\item If $n$ is even,
let $\rho\in\cale(G)_s$ such that $C_{G^*}(s)=\GL_1(q^n)$.

\item If $n$ is odd and $n\geq 3$,
let $\rho\in\cale(G)_s$ such that $C_{G^*}(s)=\rmU_1(q^n)$.
\end{enumerate}
For each case, such a semisimple element $s$ exists, and by Corollary~\ref{0314},
we have $\Theta\text{\rm -rk}(\rho)=n$.
\end{enumerate}
\end{proof}

\subsection{$\Theta$-ranks and parabolic induction}
Let $\bfG_n$ be one of $\Sp_{2n}$, $\rmO^+_{2n}$, $\rmO^-_{2n+2}$, $\rmO_{2n+1}$,
$\rmU_{2n}$, or $\rmU_{2n+1}$.
For $\rho\in\cale(G_n)$ and an integer $m\geq 0$,
recall that a subset $\cale_{G_n}^{G_m}(\rho)\subset\cale(G_m)$ is defined in Subsection~\ref{0211}.

\begin{lem}\label{0315}
Let $\rho\in\cale(G_n)_1$ and $m\geq n$.
Then
\[
\Theta\text{\rm -rk}(\rho)
=\min\{\,\Theta\text{\rm -rk}(\rho')\mid \rho'\in\cale_{G_n}^{G_m}(\rho)\,\}
\]
\end{lem}
\begin{proof}
Let $\rho$ be a unipotent character of $G_n$.
It is clear that we need only to show that
$\Theta\text{\rm -rk}(\rho')\geq\Theta\text{\rm -rk}(\rho)$ for any $\rho'\in\cale_{G_n}^{G_{n+1}}(\rho)$,
and there exists $\rho_0\in\cale_{G_n}^{G_{n+1}}(\rho)$ such that
$\Theta\text{\rm -rk}(\rho_0)=\Theta\text{\rm -rk}(\rho)$.

Now we can write $\rho=\rho_\Lambda$ ($\rho=\rho_\Lambda$ or $\rho=\rho_\Lambda\cdot\sgn$ if $\bfG=\rmO_{2n+1}$)
for some symbol $\Lambda=\binom{a_1,a_2,\ldots,a_{m_1}}{b_1,b_2,\ldots,b_{m_2}}\in\cals_{\bfG_n}$.
Consider the following reductive dual pairs $(\bfG'_{n'},\bfG_n)$:
\begin{enumerate}
\item[(a)] $(\rmO^+_{2n'},\Sp_{2n})$, $(\Sp_{2n'},\rmO^+_{2n})$,
$(\rmU_{2n'},\rmU_{2n})$, $(\rmU_{2n'+1},\rmU_{2n+1})$;

\item[(b)] $(\rmO^-_{2n'},\Sp_{2n})$, $(\Sp_{2n'},\rmO^-_{2n+2})$,
$(\rmU_{2n'+1},\rmU_{2n})$, $(\rmU_{2n'},\rmU_{2n+1})$.
\end{enumerate}
Let $\rho'\in\cale_{G_n}^{G_{n+1}}(\rho)$.
Then $\rho'$ is also unipotent and we write $\rho'=\rho_{\Lambda'}$
($\rho'=\rho_{\Lambda'}$ or $\rho'=\rho_{\Lambda'}\cdot\sgn$ if $\bfG_n=\rmO_{2n+1}$) where
$\Lambda'$ is equal to one of the following:
\[
\textstyle
\binom{a_1,\ldots,a_{i-1},a_i+1,a_{i+1},\ldots,a_{m_1}}{b_1,b_2,\ldots,b_{m_2}},
\binom{a_1+1,a_2+1,\ldots,a_{m_1}+1,1}{b_1+1,b_2+1,\ldots,b_{m_2}+1,0},
\binom{a_1,a_2,\ldots,a_{m_1}}{b_1,\ldots,b_{i-1},b_i+1,b_{i+1},\ldots,b_{m_2}},
\binom{a_1+1,a_2+1,\ldots,a_{m_1}+1,0}{b_1+1,b_2+1,\ldots,b_{m_2}+1,1}
\]
for some index $i$.
By Proposition~\ref{0306}, it is clear that
$\Theta\text{\rm -rk}(\rho_{\Lambda'})
=\Theta\text{\rm -rk}(\Lambda')
\geq\Theta\text{\rm -rk}(\Lambda)
=\Theta\text{\rm -rk}(\rho_\Lambda)$.
Now define $\Lambda_0\in\cals_{\bfG_{n+1}}$ by
\[
\Lambda_0=\begin{cases}
\binom{a_1+1,a_2,\ldots,a_{m_1}}{b_1,b_2,\ldots,b_{m_2}}, & \text{if $(\bfG'_{n'},\bfG_n)$ is in case (a)};\\
\binom{a_1,a_2,\ldots,a_{m_1}}{b_1+1,b_2,\ldots,b_{m_2}}, & \text{if $(\bfG'_{n'},\bfG_n)$ is in case (b)}.
\end{cases}
\]
Then $\rho_0:=\rho_{\Lambda_0}\in\cale_{G_n}^{G_{n+1}}(\rho)$ and by Proposition~\ref{0306}, we know that
$\Theta\text{\rm -rk}(\rho_{\Lambda_0})=\Theta\text{\rm -rk}(\rho_\Lambda)$.
\end{proof}

From the proof of the lemma, we have the following result:

\begin{cor}
Let $\rho\in\cale(G_n)$ and $\rho'\in\cale_{G_n}^{G_{n+1}}(\rho)$.
Then either $\Theta\text{\rm -rk}(\rho')=\Theta\text{\rm -rk}(\rho)$
or $\Theta\text{\rm -rk}(\rho')=\Theta\text{\rm -rk}(\rho)+2$.
\end{cor}

Then we generalize Lemma~\ref{0315} to any irreducible characters.

\begin{lem}\label{0301}
Let $\rho\in\cale(G_n)$ and $m\geq n$.
Then
\[
\Theta\text{\rm -rk}(\rho)
=\min\{\,\Theta\text{\rm -rk}(\rho')\mid \rho'\in\cale_{G_n}^{G_m}(\rho\chi),\ \chi\text{ a linear character of }G_n\,\}
\]
\end{lem}
\begin{proof}
As in the proof of Lemma~\ref{0315}, we need only to show that
$\Theta\text{\rm -rk}(\rho')\geq\Theta\text{\rm -rk}(\rho)$ for any $\rho'\in\cale_{G_n}^{G_{n+1}}(\rho\chi)$
for any linear character $\chi$,
and there exists $\rho_0\in\cale_{G_n}^{G_{n+1}}(\rho\chi)$ for some linear character $\chi$ such that
$\Theta\text{\rm -rk}(\rho_0)=\Theta\text{\rm -rk}(\rho)$.

Now we suppose that $\rho\in\cale(G_n)_s$ for some $s$.
When $\bfG_n$ is a symplectic group or an orthogonal group,
we write $\Xi_s(\rho)=\rho^{(0)}\otimes\rho_{\Lambda^{(-)}}\otimes\rho_{\Lambda^{(+)}}$
for some $\Lambda^{(\varepsilon)}\in\cals_{\bfG^{(\varepsilon)}}$ where $\varepsilon=\pm$.

\begin{enumerate}
\item Suppose that $\bfG_n=\Sp_{2n}$ and $\Theta\text{\rm -rk}(\rho)$ is even.
By Proposition~\ref{0307}, we see that $\Theta\text{\rm -rk}(\rho)=2n-2n^{(+)}+\Theta(\rho_{\Lambda^{(+)}})$.
Let $\rho'\in\cale_{G_n}^{G_{n+1}}(\rho)$.
Then we know that $\rho'\in\cale(G_{n+1})_{s_{n+1}}$ where
$s_{n+1}=(s,1)\in(G_n^*)^0\times T_1\subset (G^*_{n+1})^0$.
Now $\Xi_{s_{n+1}}(\rho')=\rho^{(0)}\otimes\rho_{\Lambda^{(-)}}\otimes\rho_{\Lambda'^{(+)}}$ where
$\rho_{\Lambda'^{(+)}}\in\cale^{G^{(+)}_{n^{(+)}+1}}_{G^{(+)}_{n^{(+)}}}(\rho_{\Lambda^{(+)}})$ by Lemma~\ref{0202}.
Then by Lemma~\ref{0315},
\begin{align*}
\Theta\text{\rm -rk}(\rho')
= 2(n+1)-2(n^{(+)}+1)+\Theta\text{\rm -rk}(\rho_{\Lambda'^{(+)}})
&\geq 2n-2n^{(+)}+\Theta\text{\rm -rk}(\rho_{\Lambda^{(+)}}) \\
&= \Theta\text{\rm -rk}(\rho).
\end{align*}
Next let $\Lambda^{(+)}_0$ be defined as $\Lambda_0$ in the proof of Lemma~\ref{0315},
and let
\[
\rho_0=\Xi^{-1}_{s_{n+1}}(\rho^{(0)}\otimes\rho_{\Lambda^{(-)}}\otimes\rho_{\Lambda^{(+)}_0}).
\]
Then we have $\rho_0\in\cale_{G_n}^{G_{n+1}}(\rho)$ and $\Theta\text{\rm -rk}(\rho)=\Theta\text{\rm -rk}(\rho_0)$

\item Suppose that $\bfG_n=\Sp_{2n}$ and $\Theta\text{\rm -rk}(\rho)$ is odd.
Then by Proposition~\ref{0307}, we have
$\Theta\text{\rm -rk}(\rho)=2n-2n^{(-)}+\Theta(\rho_{\Lambda^{(-)}})+1$.
The remaining proof is similar to case (1) with only ``$+$'' replaced by ``$-$''.

\item Suppose that $\bfG_n$ is an orthogonal group.
By twisting a linear character $\chi$ of $\bfG_n$,
by Lemma~\ref{0207} and Proposition~\ref{0307}, we may assume $\Theta\text{\rm -rk}(\rho)=2n-2n^{(+)}+\Theta(\rho_{\Lambda^{(+)}})$.
Then the remaining proof is similar to case (1).
\end{enumerate}

Now suppose that $\bfG_n$ is a unitary group.
Write $\grL_s(\rho)=\bigotimes_{i=1}^r\bigotimes_{j=1}^{t_i}\rho^{(ij)}$ as in (\ref{0308}).
By Proposition~\ref{0311},
after twisting a linear character $\chi$ of $\bfG_n$ if necessary,
we may assume that $\Theta\text{\rm -rk}(\rho)=n-n_{1,1}+\Theta(\rho^{(1,1)})$
where the notations are as in Subsection~\ref{0316}.
Then the remaining proof is similar to case (1), again.
\end{proof}

\section{$\Theta$-ranks and Other Ranks}

\subsection{$U$-rank, $A$-rank, $\overline A$-rank, and $\otimes$-rank}\label{0414}
First we recall the definition of the $U$-rank and the asymptotic rank of an irreducible character $\rho$
of a classical group $G=G(V)$ given in \cite{gurevich-howe} and \cite{gurevich-howe-rank}
where $V$ is the standard space with an associated form over $\bff_q$
(if $G$ is a symplectic group or an orthogonal group) or $\bff_{q^2}$ (if $G$ is a unitary group)
on which $G$ acts.
Suppose that we have a Witt decomposition
\[
V=X\oplus V_0\oplus Y
\]
where $X,Y$ are totally isotropic and dual to each other, and $V_0$ is anisotropic.
The group
\[
U=\{\,g\in G(V)\mid(g-1)(X\oplus V_0)=0,\ (g-1)(Y)\subset X\,\}
\]
is an abelian subgroup of $G$.
Then we have a decomposition
\[
\rho|_U=\sum_{\phi\in\cale(U)}m_\phi\phi
\]
for $m_\phi\in\bbN\cup\{0\}$.
We say that $\rho|_U$ \emph{contains} $\phi$ if $m_\phi\neq 0$.
Now each character $\phi$ of $U$ is naturally assigned a rank.

The irreducible character $\rho$ is said to be of \emph{$U$-rank} $k$, denoted by $U\text{\rm -rk}(\rho)=k$,
if the restriction $\rho|_U$ contains a character $\phi$ of rank $k$,
but contains no character of higher rank.
From the definition, we see that $0\leq U\text{\rm -rk}(\rho)\leq n$ where $n=\dim(X)=\dim(Y)$.
Moreover, it is clear that $U\text{\rm -rk}(\rho)=U\text{\rm -rk}(\rho\chi)$
for any linear character $\chi$ of $G$ because $\chi|_U$ is always trivial.
An irreducible character $\rho\in\cale(G_n)$ is said to be of \emph{low $U$-rank} if
\[
U\text{\rm -rk}(\rho)<\begin{cases}
n, & \text{if $\bfG_n=\Sp_{2n}$, $\rmU_{2n}$ or $\rmU_{2n+1}$};\\
n-1, & \text{if $\bfG_n=\rmO^+_{2n}$, $\rmO^-_{2n+2}$, or $\rmO_{2n+1}$}.
\end{cases}
\]

An irreducible character $\rho$ of $G_n$ is said to have \emph{asymptotic rank} less than or equal to $k$,
if for all sufficiently large $N$, there is a character $\rho_N\in\cale(G_N)$
such that $U\text{\rm -rk}(\rho_N)\leq k<N$ and $\rho$ occurs in $\rho_N|_{G_n}$.
The asymptotic rank of $\rho$ is denoted by $A\text{\rm -rk}(\rho)$.
As indicated in \cite{gurevich-howe-rank}, it is known that
\begin{equation}\label{0410}
U\text{\rm -rk}(\rho)\leq A\text{\rm -rk}(\rho)
\end{equation}
for any $\rho\in\cale(G)$.

To study the representation theory of reductive groups,
parabolic induction is usually more convenient than usual induction,
so we modified the definition of asymptotic rank to define the parabolic asymptotic rank as follows.
An irreducible character $\rho$ of $G_n$ is said to have \emph{parabolic asymptotic rank} less than or equal to $k$,
if for all sufficiently large $N$, there is a character $\rho_N\in\cale_{G_n}^{G_N}(\rho)$
such that $U\text{\rm -rk}(\rho_N)\leq k<N$.
The parabolic asymptotic rank of $\rho$ is denoted by $\overline A\text{\rm -rk}(\rho)$.
It is clear from the definition that
\begin{equation}\label{0411}
A\text{\rm -rk}(\rho)\leq\overline A\text{\rm -rk}(\rho)
\end{equation}
for any $\rho\in\cale(G)$.

\begin{lem}
Let $\rho$ be an irreducible character of $G_n$ and $m\geq n$.
Then
\[
\overline A\text{\rm -rk}(\rho)=\min\{\,\overline A\text{\rm -rk}(\rho')\mid \rho'\in\cale_{G_n}^{G_m}(\rho)\,\}.
\]
\end{lem}
\begin{proof}
Let $\rho'\in\cale_{G_n}^{G_m}(\rho)$ and suppose that $\overline A\text{\rm -rk}(\rho')=k_1$.
This means for all sufficiently large $N$, there exists $\rho_N\in\cale^{G_N}_{G_m}(\rho')$
such that $U\text{\rm -rk}(\rho_N)=k_1<N$.
Clearly, $\rho_N$ is in $\cale_{G_n}^{G_N}(\rho)$,
so we have $\overline A\text{\rm -rk}(\rho)\leq k_1$.

On the other hand, suppose that $\overline A\text{\rm -rk}(\rho)=k_2$.
Then for all sufficiently large $N$, there exists $\rho_N\in\cale^{G_N}_{G_n}(\rho)$
such that $U\text{\rm -rk}(\rho_N)=k_2<N$.
For $m\geq n$, we know that
\[
\cale^{G_N}_{G_n}(\rho)=\bigcup_{\rho'\in\cale^{G_m}_{G_n}(\rho)}\cale^{G_N}_{G_m}(\rho'),
\]
i.e., $\rho_N\in\cale^{G_N}_{G_m}(\rho')$ for some $\rho'\in\cale^{G_m}_{G_n}(\rho)$.
Then we have $\overline A\text{\rm -rk}(\rho')\leq k_2$ for some $\rho'\in\cale^{G_m}_{G_n}(\rho)$.
Thus the lemma is proved.
\end{proof}

According to \cite{gurevich-howe-rank} definition 6.4.1,
the \emph{tensor rank} of an irreducible character $\rho\in\cale(G)$, denoted by $\otimes\text{\rm -rk}(\rho)$,
is defined as follows:
\begin{enumerate}
\item Suppose that $\bfG$ is a symplectic group.
An irreducible character $\rho\in\cale(G)$ has tensor rank less than or equal to $k$ if
$\rho$ occurs in
\[
(\omega_{\rmO^{\epsilon_1}_1,\bfG}\otimes\omega_{\rmO^{\epsilon_2}_1,\bfG}\otimes
\cdots\otimes\omega_{\rmO^{\epsilon_k}_1,\bfG})|_G
\]
for some $\epsilon_i=+$ or $-$ and $i=1,\ldots,k$.

\item Suppose that $\bfG$ is an orthogonal group.
An irreducible character $\rho\in\cale(G)$ has tensor rank less than or equal to $2k$ if
$\rho\chi$ occurs in $(\omega_{\Sp_2,\bfG})^{\otimes k}|_G$ for some linear character $\chi$ of $G$.

\item Suppose that $\bfG$ is a unitary group.
An irreducible character $\rho\in\cale(G)$ has tensor rank less than or equal to $k$ if
$\rho\chi$ occurs in $(\omega_{\rmU_1,\bfG})^{\otimes k}|_G$ for some linear character $\chi$ of $G$.
\end{enumerate}
Here $\omega^{\otimes k}$ means the tensor product $\omega\otimes\cdots\otimes\omega$
of $k$ copies of $\omega$.
Because
\begin{align*}
(\omega_{\rmO^{\epsilon_1}_1,\bfG}\otimes\omega_{\rmO^{\epsilon_2}_1,\bfG}\otimes
\cdots\otimes\omega_{\rmO^{\epsilon_k}_1,\bfG})|_G
&= (\omega_{\rmO^\epsilon_k,\bfG})|_G \\
(\omega_{\Sp_2,\bfG})^{\otimes k}|_G
&= (\omega_{\Sp_{2k},\bfG})|_G \\
(\omega_{\rmU_1,\bfG})^{\otimes k}|_G
&=(\omega_{\rmU_k,\bfG})|_G
\end{align*}
for some $\epsilon=+$ or $-$,
we see that
\begin{equation}\label{0412}
\otimes\text{\rm -rk}(\rho)=\Theta\text{\rm -rk}(\rho)
\end{equation}
for any $\rho\in\cale(G)$.

\subsection{$\eta$- and $\underline\theta$-correspondence}
The following fundamental result is from \cite{gurevich-howe-rank}:

\begin{prop}[Gurevich-Howe]\label{0401}
Let $(\bfG',\bfG)$ be a dual pair in stable range, i.e., it is of one of the following types:
\begin{enumerate}
\item[(i)] $(\Sp_k,\rmO^+_{2n})$, $(\Sp_k,\rmO^-_{2n+2})$ or $(\Sp_k,\rmO_{2n+1})$ with $k$ even and $k\leq n$;

\item[(ii)] $(\rmO_k^\epsilon,\Sp_{2n})$ with $k\leq n$;

\item[(iii)] $(\rmU_k,\rmU_{2n})$ or $(\rmU_k,\rmU_{2n+1})$ with $k\leq n$.
\end{enumerate}
Then
for $\rho'\in\cale(G')$ there is a unique irreducible character
$\eta(\rho')\in\Theta_\bfG(\rho')$ of $U$-rank $k$, and
all other elements in $\Theta_\bfG(\rho')$ have $U$-rank less than $k$.
Furthermore, the mapping $\eta\colon\cale(G')\rightarrow\cale(G)$ is injective.
\end{prop}

The mapping $\eta\colon\cale(G')\rightarrow\cale(G)$
is called the \emph{$\eta$-correspondence} for the dual pair $(\bfG',\bfG)$
in stable range.
In \cite{pan-eta} and \cite{pan-eta-unitary}, a sub-relation of $\Theta$,
called the \emph{$\underline\theta$-correspondence}, is defined as follows.
Let $(\bfG',\bfG)$ be one of the following types of dual pairs:
\begin{enumerate}
\item[(I)] $(\rmO^+_{2n'},\Sp_{2n})$ or $(\Sp_{2n'},\rmO^+_{2n})$

\item[(II)] $(\rmO^-_{2n'},\Sp_{2n})$ or $(\Sp_{2n'},\rmO^-_{2n})$

\item[(III)] $(\rmU_{2n'},\rmU_{2n})$ or $(\rmU_{2n'+1},\rmU_{2n+1})$

\item[(IV)] $(\rmU_{2n'},\rmU_{2n+1})$ or $(\rmU_{2n'+1},\rmU_{2n})$
\end{enumerate}
Let $\Lambda'\in\cals_{\bfG'}$ and $\Lambda\in\cals_{\bfG}$
such that $(\Lambda',\Lambda)\in\calb_{\bfG',\bfG}$.
Then we know that $|\Upsilon(\Lambda)|=|\Upsilon(\Lambda')|+\tau$
where $\tau$ is given by
\[
\tau=\begin{cases}
n-n'+\left\lfloor\frac{{\rm def}(\Lambda)}{2}\right\rfloor, & \text{for case (I);}\\
n-n'-\left\lfloor\frac{{\rm def}(\Lambda)}{2}\right\rfloor, & \text{for case (II);}\\
n-n'+|{\rm def}(\Lambda)|, & \text{for case (III);}\\
n-n'-|{\rm def}(\Lambda)|-1, & \text{for case (IV).}\\
\end{cases}
\]
For $\Lambda'\in\cals_{\bfG'}$ such that $\tau\geq 0$,
we define $\underline\theta(\Upsilon(\Lambda'))$ by
\[
\underline\theta\left(\sqbinom{\mu_1,\mu_2,\ldots,\mu_{m_1}}{\nu_1,\nu_2,\ldots,\nu_{m_2}}\right)
=\begin{cases}
\sqbinom{\nu_1,\nu_2,\ldots,\nu_{m_2}}{\mu_1,\mu_2,\ldots,\mu_{m_1}}\cup\sqbinom{\tau}{-},
& \text{for cases (I),(III) ;}\\
\sqbinom{\nu_1,\nu_2,\ldots,\nu_{m_2}}{\mu_1,\mu_2,\ldots,\mu_{m_1}}\cup\sqbinom{-}{\tau},
& \text{for cases (II),(IV).}
\end{cases}
\]
Then we define $\underline\theta(\Lambda')$ to be the symbol $\Lambda\in\cals_\bfG$ such that
$\Upsilon(\Lambda)=\underline\theta(\Upsilon(\Lambda'))$.
Hence we obtain a relation $\underline\theta$ between $\cale(G')_1$ and $\cale(G)_1$
by letting $\underline\theta(\rho_{\Lambda'})=\rho_{\underline\theta(\Lambda')}$.

Next we extend the domain of $\underline\theta$ by letting it compatible with the Lusztig
correspondence, i.e.,
\begin{enumerate}
\item for $(\bfG',\bfG)$ being of one cases (I),(II),(III),(IV) above, $\rho'\in\cale(G')_{s'}$
and $\rho\in\cale(G)_s$, we have the following commutative diagram:
\[
\begin{CD}
\rho' @> \underline\theta_{\bfG',\bfG} >> \rho \\
@V \Xi_{s'} VV @VV \Xi_s V \\
\rho'^{(0)}\otimes\rho_{\Lambda'^{(-)}}\otimes\rho_{\Lambda'^{(+)}}
@> \id\otimes\id\otimes\underline\theta_{\bfG'^{(+)},\bfG^{(+)}} >>
\rho^{(0)}\otimes\rho_{\Lambda^{(-)}}\otimes\rho_{\Lambda^{(+)}}
\end{CD}
\]

\item for the dual pair $(\bfG',\bfG)=(\rmO_{2n'+1},\Sp_{2n})$, $\rho'\in\cale(G')_{s'}$
and $\rho\in\cale(G)_s$, we have the following commutative diagram
\[
\begin{CD}
\rho' @> \underline\theta_{\bfG',\bfG} >> \rho \\
@V \iota\circ\Xi_{s'} VV @VV \Xi_s V \\
\rho'^{(0)}\otimes\rho_{\Lambda'^{(+)}}\otimes\rho_{\Lambda'^{(-)}}
@> \id\otimes\underline\theta_{\bfG'^{(+)},\bfG^{(-)}}\otimes\id >>
\rho^{(0)}\otimes\rho_{\Lambda^{(-)}}\otimes\rho_{\Lambda^{(+)}}
\end{CD}
\]
\end{enumerate}
We know that $\underline\theta$ is a sub-correspondence of $\Theta$,
and we can also define the \emph{$\underline\theta$-rank} of $\rho\in\cale(G)$ by
the same way as $\Theta$-rank.
It is know that 
\begin{equation}\label{0415}
\underline\theta\text{\rm -rk}(\rho')=\Theta\text{\rm -rk}(\rho')
\end{equation}
for any $\rho\in\cale(G')$ by \cite{pan-eta} lemma~6.11.

\begin{prop}\label{0407}
Let $(\bfG',\bfG)$ be a reductive dual pair in stable range.
Then $\eta=\underline\theta$.
\end{prop}
\begin{proof}
If $(\bfG',\bfG)$ is a symplectic/orthogonal dual pair in stable range,
this proposition is \cite{pan-eta} theorem 6.16.
Moreover, it is not difficult to see that the proof in \cite{pan-eta} works for a unitary dual pair in table range.
\end{proof}

\begin{lem}\label{0504}
Let $\rho\in\cale(G_n)$.
If\/ $\Theta\text{\rm -rk}(\rho)\leq n$,
then $U\text{\rm -rk}(\rho)=\Theta\text{\rm -rk}(\rho)$.
\end{lem}
\begin{proof}
Suppose that $\Theta\text{\rm -rk}(\rho)=k\leq n$.
So by (\ref{0415}) we have $\underline\theta\text{\rm -rk}(\rho)=k\leq n$, i.e.,
$\rho=\underline\theta(\rho')$ for some $\rho'\in\cale(\rmO^\epsilon_k(q))$ for some $\epsilon=+$ or $-$.
Now the dual pair $(\rmO^\epsilon_k,\Sp_{2n})$ is in stable range,
so the two mappings $\eta$ and $\underline\theta$ coincides by Proposition~\ref{0407},
i.e., $\rho=\eta(\rho')$ and hence $U\text{\rm -rk}(\rho)=k$ by Proposition~\ref{0401}.
\end{proof}

\begin{cor}\label{0408}
We have $\overline A\text{\rm -rk}(\rho)\leq \Theta\text{\rm -rk}(\rho)$ for any $\rho\in\cale(G)$.
\end{cor}
\begin{proof}
Let $\bfG=\bfG_n$, $\rho\in\cale(G)$ and suppose that $\Theta\text{\rm -rk}(\rho)=k$.
So now $\rho'\otimes\rho$ occurs in the $\Theta$-correspondence for some dual pair $(G',G_n)$.
Then we know that, for any $N\geq n$, $\rho'\otimes\rho_N$ occurs in the $\Theta$-correspondence
for the dual pair $(G',G_N)$ for some $\rho_N\in\cale_{G_n}^{G_N}(\rho)$.
By definition, we have $\Theta\text{\rm -rk}(\rho_N)\leq k$,
and by Lemma~\ref{0301}, we have $\Theta\text{\rm -rk}(\rho_N)\geq\Theta\text{\rm -rk}(\rho)=k$.
Therefore, we conclude that $\Theta\text{\rm -rk}(\rho_N)=k$.
When $N\geq k$, by Lemma~\ref{0504}, we see that $U\text{\rm -rk}(\rho_N)=k$.
Then by the definition of $\overline A$-rank,
we conclude that $\overline A\text{\rm -rk}(\rho)\leq U\text{\rm -rk}(\rho_N)=k$.
Therefore, $\overline A\text{\rm -rk}(\rho)\leq \Theta\text{\rm -rk}(\rho)$.
\end{proof}

\begin{proof}[Proof of Theorem~\ref{0101}]
Now from (\ref{0410}), (\ref{0411}), (\ref{0412}) and Corollary~\ref{0408},
we can summarize the relation between these various ranks as
\begin{equation}\label{0413}
U\text{\rm -rk}(\rho)\leq A\text{\rm -rk}(\rho)\leq \overline A\text{\rm -rk}(\rho)\leq \otimes\text{\rm -rk}(\rho)
=\Theta\text{\rm -rk}(\rho)
\end{equation}
for any $\rho\in\cale(G)$.
Then Theorem~\ref{0101} follows from Lemma~\ref{0504} immediately.
\end{proof}

\subsection{Agreement of $\overline A$-ranks and $\Theta$-ranks}
It is proved by Howe-Gurevich in \cite{gurevich-howe-rank} theorem 11.4
that $A\text{\rm -rk}(\rho)=\otimes\text{\rm -rk}(\rho)$
for any $\rho\in\cale(\GL_n(q))$.
Now we want to apply their method to other classical groups.

Let $G=G_n=G(V)$ and $V=X_n\oplus V_0\oplus Y_n$ be a complete polarization of $V$ where $V_0$ is anisotropic.
Let $\{x_1,x_2,\ldots,x_n\}$ be a basis of $X_n$, and let $\{y_1,y_2,\ldots,y_n\}$ be the
corresponding dual basis of $Y_n$, i.e., $\langle x_i,y_j\rangle=\delta_{ij}$
where $\delta_{ij}$ is the Kronecker delta symbol.
For $1\leq k\leq n$, let $X_k$ (resp.~$Y_k$) denote the space spanned by $\{x_1,x_2,\ldots,x_k\}$
(resp.~$\{y_1,y_2,\ldots,y_k\}$).
So we have the decomposition $V=X_k\oplus V_{n-k}\oplus Y_k$ where $V_{n-k}$ is a space of Witt index $n-k$.
Define
\begin{align*}
H_k &=\{\,g\in G(V)\mid(g-1)(X_k)=0,\ (g-1)(Y_k)=0\,\}.
\end{align*}
Then $H_k$ is a subgroup of $G(V)$, $H_k\subset H_{k-1}$ for each $k$, and $H_0=G_n$.

\begin{lem}\label{0406}
Let $\rho\in\cale(G_n)$ be of low $U$-rank $k$, i.e.,
\[
U\text{\rm -rk}(\rho)=k<\begin{cases}
n, & \text{if\/ $\bfG_n$ is $\Sp_{2n}$, $\rmU_{2n}$ or $\rmU_{2n+1}$};\\
n-1, & \text{if\/ $\bfG_n$ is $\rmO^+_{2n}$, $\rmO^-_{2n+2}$ or $\rmO_{2n+1}$}.
\end{cases}
\]
Then $\rho|_{H_k}$ contains a linear character.
\end{lem}
\begin{proof}
Suppose that $\rho\in\cale(G_n)$ is of low $U$-rank $k$.
We define the following unipotent subgroups of $G(V)$:
\begin{align*}
U_1 &=\{\,g\in G(V)\mid(g-1)(X_k)=0,\ (g-1)(V_{n-k}\oplus Y_k)\subset X_k \,\}, \\
U_2 &=\{\,g\in G(V)\mid(g-1)(X_k\oplus V_{n-k})=0,\ (g-1)(Y_k)\subset X_k\,\}, \\
U_3 &=U_1\cap H_k.
\end{align*}
It is not difficult to see that $U_1\simeq U_2\times U_3$,
and $U_1$ is normalized by $H_k$.
Define the subgroup $U_4$ of $G(V)$ as follows:
\begin{itemize}
\item If $G(V)$ is a symplectic group,
then $k<n$,
and let $U_4$ be the subgroup of elements $g\in G(V)$ such that
\begin{itemize}
\item $(g-1)y_{k+1}=a x_{k+1}$, for $a\in\bff_q$, and

\item $(g-1)(X_{k+1}\oplus V_{n-k-1}\oplus Y_k)=0$.
\end{itemize}

\item If $G(V)$ is an orthogonal group,
then $k<n-1$,
and let $U_4$ be the subgroup of elements $g\in G(V)$ such that
\begin{itemize}
\item $(g-1)y_{k+2}=a x_{k+1}$ and $(g-1)y_{k+1}=-a x_{k+2}$, for $a\in\bff_q$, and

\item $(g-1)(X_{k+2}\oplus V_{n-k-2}\oplus Y_k)=0$.
\end{itemize}

\item If $G(V)$ is a unitary group,
then $k<n$,
and let $U_4$ be the subgroup of elements $g\in G(V)$ such that
\begin{itemize}
\item $(g-1)y_{k+1}=a x_{k+1}$, for $a\in\bff_{q^2}$ satisfying $\overline a=-a$, and

\item $(g-1)(X_{k+1}\oplus V_{n-k-1}\oplus Y_k)=0$.
\end{itemize}
Here $a\mapsto\overline a$ is the nontrivial automorphism of $\bff_{q^2}$ over $\bff_q$.
\end{itemize}
Finally we define the subgroup $U_5$ of $G(V)$ as follows:
\begin{itemize}
\item If $G(V)$ is a symplectic group or a unitary group,
let $U_5$ be the subgroup of elements $g\in G(V)$ such that
\begin{itemize}
\item $(g-1)(X_{k+1}\oplus V_{n-k-1})=0$ and

\item $(g-1)(Y_{k+1})\subset X_{k+1}$.
\end{itemize}

\item If $G(V)$ is an orthogonal group,
let $U_5$ be the subgroup of elements $g\in G(V)$ such that
\begin{itemize}
\item $(g-1)(X_{k+2}\oplus V_{n-k-2})=0$ and

\item $(g-1)(Y_{k+2})\subset X_{k+2}$.
\end{itemize}
\end{itemize}

There exists a linear character $\psi$ of $U_1$ of rank $k$ and trivial on $U_3$,
so $\psi|_{U_2}$ still has rank $k$.
It is easy to see that $h(u_2,u_3)h^{-1}=(u_2,u_3')$ for $h\in H_k$, $u_2\in U_2$ and $u_3\in U_3$.
Because now $\psi$ is trivial on $U_3$, $\psi$ is invariant under the conjugation by $H_k$.
Because now $U\text{\rm -rk}(\rho)=k$,
$\rho|_{U_1}$ contains a nontrivial $\psi$-eigenspace, denoted by $\rho_\psi$.
Because $\psi$ is invariant under the conjugation of $H_k$, the space
$\rho_\psi$ is invariant under the action of $\rho|_{H_k}$.
Let $\sigma$ denote the representation of $H_k$ on the space $\rho_\psi$.

Now we claim that the representation $\sigma$ of $H_k$ is trivial on the commutator subgroup
of $H_k$.
Suppose that the claim is not true, i.e.,
$\sigma$ is not trivial on the commutator subgroup of $H_k$
Then $\sigma$ is not trivial on $U_4$,
because the normal subgroup of $H_k$ generated by $U_4$ is the full commutator subgroup of $H_k$.
The space $\rho_\psi$ is decomposed into eigenspaces for $U_4$.
Let $\rho_\psi^1$ be the eigenspace corresponding to a non-trivial linear character $\lambda$ of $U_4$.
The group $U_5$ is commutative and contains both $U_2$ and $U_4$.
So the space $\rho_\psi^1$ is decomposed into eigenspaces of $U_5$.
The character of $U_5$ which corresponding a nontrivial eigenspace in $\rho_\psi^1$ of $U_5$
has rank greater than $k$ and we get a contradiction.

Now the representation $\sigma$ of $H_k$ is trivial on the commutator subgroup
of $H_k$, we see that $\sigma$ is a sum of linear characters of $H_k$.
\end{proof}

Recall that under the \emph{Shr\" odinger model} of the Weil representation $\omega$
(\cf.~\cite{gurevich-howe}),
the group $G(V)$ acts on the space $L^2(\Hom(X_k,V))$ of functions on $\Hom(X_k,V )$ by
\[
(\omega(g)f)(T)=\chi_k(g) f(g^{-1}\circ T)
\]
for $g\in G(V)$ and $T\in\Hom(X_k,V)$ where $\chi_k$ is a linear character of
$G(V)$ of order at most $2$.
Recall that $L^2(G/H_k,\chi_k^{-1})$ denotes the space of functions $f$ on $G$ such that
$f(gh)=\chi_k^{-1}(h)f(g)$ for $g\in G$ and $h\in H_k$.
Then $G(V)$ acts on the space $L^2(G/H_k,\chi_k)$ by
\[
(\pi(g)f)(x)=f(g^{-1}x)
\]
for $g,x\in G(V)$.

Fix an element $T_0\in\Hom(X_k,V)$ such that the image of $T_0$ is equal to $X_k$.
Note that for $h\in H_k$, we have $h\circ T_0=T_0$ from the definition of $H_k$.
Then it is not difficult to see that the mapping
\[
\phi\colon L^2(\Hom(X_k,V))\rightarrow L^2(G/H_k,\chi_k)
\]
given by $\phi(f)(x)=\chi_k(x^{-1})f(x\circ T_0)$ is a surjective linear transformation.
Moreover, we have
\begin{align*}
\pi(g')(\phi(f))(g) &=\phi(f)(g'^{-1}g)=\chi_k(g^{-1}g')f(g'^{-1}g\circ T) \\
\phi(\omega(g')f)(g) &=\chi_k(g^{-1})(\omega(g')f)(g\circ T)=\chi_k(g^{-1}g')f(g'^{-1}g\circ T),
\end{align*}
i.e., $\phi$ commutes the two actions $\omega$ and $\pi$ of $G$.
Therefore $\pi$ can be regarded as a subrepresentation of $\omega$.

\begin{lem}\label{0409}
Let $\rho\in\cale(G_n)$ and $0\leq k\leq n$.
If $\rho|_{H_k}$ contains a linear character,
then $\Theta\text{\rm -rk}(\rho)\leq 2k$.
\end{lem}
\begin{proof}
Suppose that $\rho\in\cale(G_n)$ and $\rho|_{H_k}$ contains a linear character.
Since a linear character of $H_k$ is the restriction of a linear character of $G_n$,
after twisting a linear character of $G_n$, we may assume that
$\rho|_{H_k}$ contains the character $\chi_k^{-1}$.
By Frobenius reciprocity, $\rho$ is contained in $L^2(G_n/H_k,\chi_k^{-1})$
under the action $\pi$ of left translation of $G_n$.
Then by the above discussion, we see that $\rho$ is contained in
the Sch\"ondinger model $L^2(\Hom(X_k,V))$ of the Weil representation $\omega$.
This means that $\rho$ occurs in the $\Theta$-correspondence for the dual pair $(\bfG',\bfG)$
where $\bfG'$ is the group acting on a space $V'$ which has a Witt decomposition with a trivial anisotropic kernel and
a maximal totally isotropic space isomorphic to $X_k$, i.e., $\dim(V')=2k$.
This concludes that $\Theta\text{\rm -rk}(\rho)\leq 2k$.
\end{proof}

Now we have an analogous result of \cite{gurevich-howe-rank} proposition 11.2.

\begin{cor}
If $\rho\in\cale(G_n)$ is of $U$-rank $k\leq\frac{n}{2}$,
then $U\text{\rm -rk}(\rho)=\Theta\text{\rm -rk}(\rho)$.
\end{cor}
\begin{proof}
Suppose that $U\text{\rm -rk}(\rho)=k\leq \frac{n}{2}$.
Then by Lemma~\ref{0406} and Lemma~\ref{0409} we have $\Theta\text{\rm -rk}(\rho)\leq 2k\leq n$.
Then by Lemma~\ref{0504}, we conclude that $\Theta\text{\rm -rk}(\rho)=U\text{\rm -rk}(\rho)=k$.
\end{proof}

\begin{proof}[Proof of Theorem~\ref{0102}]
Let $\rho\in\cale(G_n)$ where $\bfG_n$ is $\Sp_{2n}$, $\rmO^+_{2n}$, $\rmO^-_{2n+2}$, $\rmO_{2n+1}$, $\rmU_{2n}$,
or $\rmU_{2n+1}$.
Suppose that $\overline A\text{\rm -rk}(\rho)=k_1$ and $\Theta\text{\rm -rk}(\rho)=k_2$.
We know that $k_1\leq k_2$ by Corollary~\ref{0408} and $k_2\leq 2n+2$ from Remark~\ref{0304}.
Then by definition of parabolic asymptotic rank we have
\begin{equation}\label{0404}
k_1=\min\{\,U\text{\rm -rk}(\rho_N)\mid\rho_N\in\cale_{G_n}^{G_N}(\rho\chi),\ \chi\text{ linear},\ N>\!>n\,\},
\end{equation}
and by Lemma~\ref{0301} we know that
\begin{equation}\label{0405}
k_2=\min\{\,\Theta\text{\rm -rk}(\rho_N)\mid\rho_N\in\cale_{G_n}^{G_N}(\rho\chi),\ \chi\text{ linear},\ N>\!>n\,\}.
\end{equation}
Let $N$ be sufficiently larger than $n$, in particular, we assume that $N\geq 4n+4$.
Then by (\ref{0404}) there is $\rho_1\in\cale_{G_n}^{G_N}(\rho\chi_1)$ for some linear character $\chi_1$ of $G_n$
such that $U\text{\rm -rk}(\rho_1)=k_1\leq k_2\leq 2n+2\leq \frac{N}{2}$.
Then by Lemma~\ref{0409}, we have $\Theta\text{\rm -rk}(\rho_1)\leq 2k_1\leq N$.
Then by Lemma~\ref{0504}, we have $\Theta\text{\rm -rk}(\rho_1)=U\text{\rm -rk}(\rho_1)=k_1$.
Then (\ref{0405}) implies that
\[
k_2=\Theta\text{\rm -rk}(\rho)\leq\Theta\text{\rm -rk}(\rho_1)=k_1=\overline A\text{\rm -rk}(\rho)\leq k_2.
\]
Thus the theorem is proved.
\end{proof}

\bibliography{refer}
\bibliographystyle{amsalpha}

\end{document}